\documentclass[sn-mathphys,iicol]{sn-jnl}

\usepackage{lmodern}
\usepackage{standalone}
\usepackage{breqn}
\usepackage{graphicx}%
\usepackage{multirow}%
\usepackage{amsmath,amssymb,amsfonts}%
\usepackage{amsthm}%
\usepackage{mathrsfs}%
\usepackage[title]{appendix}%
\usepackage[dvipsnames]{xcolor}%
\usepackage{textcomp}%
\usepackage{manyfoot}%
\usepackage{booktabs}%
\usepackage{algorithm}%
\usepackage{algorithmicx}%
\usepackage{algpseudocode}%
\usepackage{listings}%
\usepackage{hyperref}%
\usepackage{cleveref}%
\usepackage{tikz}%
\usepackage{pgfplots}
\usepackage{xspace}
\usepackage{accents}
\usepackage{listings}
\usepackage{color}
\usepackage{mathtools}
\usepackage{subcaption}
\usepackage{booktabs}
\usepackage[final]{pdfpages}
\usepackage{nicefrac}
\usepackage{mdframed}
\usepackage{gensymb}
\usepackage{setspace}
\usepackage{changepage} 
\usepackage{cuted} 
\usepackage{widetext}
\usepackage{siunitx}

\usetikzlibrary{patterns}
\usetikzlibrary{shapes,arrows}

\newcommand\cropped[1]{%
    \centering
	\includegraphics[width=0.8\linewidth]{#1cropped.png}}

\newcommand\mydots{\hbox to 1em{.\hss.\hss.}}

\makeatletter
\pgfplotsset{
	/tikz/max node/.style={
		anchor=south
	},
	/tikz/min node/.style={
		anchor=north
	},
	mark minBelow/.style={
		point meta rel=per plot,
		visualization depends on={x \as \xvalue},
		scatter/@pre marker code/.code={%
			\ifx\pgfplotspointmeta\pgfplots@metamin
			\def\markopts{}%
			\node [min node,below] {
				\pgfmathprintnumber[fixed relative,precision=4, set thousands separator={}]{\pgfkeysvalueof{/data point/y}}
			};
			\else
			\def\markopts{mark=none}
			\fi
			\expandafter\scope\expandafter[\markopts,every node near coord/.style=green]
		},%
		scatter/@post marker code/.code={%
			\endscope
		},
		scatter,
	},
	mark minAbove/.style={
	point meta rel=per plot,
	visualization depends on={x \as \xvalue},
	scatter/@pre marker code/.code={%
		\ifx\pgfplotspointmeta\pgfplots@metamin
		\def\markopts{}%
		\node [min node,above] {
			\pgfmathprintnumber[fixed relative,precision=4, set thousands separator={}]{\pgfkeysvalueof{/data point/y}}
		};
		\else
		\def\markopts{mark=none}
		\fi
		\expandafter\scope\expandafter[\markopts,every node near coord/.style=green]
	},%
	scatter/@post marker code/.code={%
		\endscope
	},
	scatter,
	},
	mark maxAbove/.style={
		point meta rel=per plot,
		visualization depends on={x \as \xvalue},
		scatter/@pre marker code/.code={%
			\ifx\pgfplotspointmeta\pgfplots@metamax
			\def\markopts{}%
			\node [max node,above] {
				\pgfmathprintnumber[fixed relative,precision=4, set thousands separator={}]{\pgfkeysvalueof{/data point/y}}
			};
			\else
			\def\markopts{mark=none}
			\fi
			\expandafter\scope\expandafter[\markopts]
		},%
		scatter/@post marker code/.code={%
			\endscope
		},
		scatter
	},
	mark maxBelow/.style={
	point meta rel=per plot,
	visualization depends on={x \as \xvalue},
	scatter/@pre marker code/.code={%
		\ifx\pgfplotspointmeta\pgfplots@metamax
		\def\markopts{}%
		\node [max node,below] {
			\pgfmathprintnumber[fixed relative,precision=4, set thousands separator={}]{\pgfkeysvalueof{/data point/y}}
		};
		\else
		\def\markopts{mark=none}
		\fi
		\expandafter\scope\expandafter[\markopts]
	},%
	scatter/@post marker code/.code={%
		\endscope
	},	scatter
},
	mark end/.style={
	point meta rel=per plot,
	visualization depends on={x \as \xvalue},
	scatter/@pre marker code/.code={%
		
		\pgfmathtruncatemacro\usemark{
			(\coordindex==(\numcoords-1)
		}
		\ifnum\usemark=1
		\def\markopts{}%
		\node [min node,right] {
			\pgfmathprintnumber[fixed relative,precision=4, set thousands separator={}]{\pgfkeysvalueof{/data point/y}}
		};
		\else
		\def\markopts{mark=none}
		\fi
		\expandafter\scope\expandafter[\markopts,every node near coord/.style=green]
	},%
	scatter/@post marker code/.code={%
		\endscope
	},	scatter,
	},
mark endAbove/.style={
point meta rel=per plot,
visualization depends on={x \as \xvalue},
scatter/@pre marker code/.code={%
	
	\pgfmathtruncatemacro\usemark{
		(\coordindex==(\numcoords-1)
	}
	\ifnum\usemark=1
	\def\markopts{}%
	\node [min node,above] {
		\pgfmathprintnumber[fixed relative,precision=4, set thousands separator={}]{\pgfkeysvalueof{/data point/y}}
	};
	\else
	\def\markopts{mark=none}
	\fi
	\expandafter\scope\expandafter[\markopts,every node near coord/.style=green]
},%
scatter/@post marker code/.code={%
	\endscope
},	scatter,
},
mark endBelow/.style={
point meta rel=per plot,
visualization depends on={x \as \xvalue},
scatter/@pre marker code/.code={%
	
	\pgfmathtruncatemacro\usemark{
		(\coordindex==(\numcoords-1)
	}
	\ifnum\usemark=1
	\def\markopts{}%
	\node [min node,below] {
		\pgfmathprintnumber[fixed relative,precision=4, set thousands separator={}]{\pgfkeysvalueof{/data point/y}}
	};
	\else
	\def\markopts{mark=none}
	\fi
	\expandafter\scope\expandafter[\markopts,every node near coord/.style=green]
},%
scatter/@post marker code/.code={%
	\endscope
},	scatter,
}
}
\makeatother

\DeclareMathOperator{\spn}{span}

\newcommand{\freefem}{\textit{FreeFEM++}\xspace}
\newcommand{\python}{\textit{Python}\xspace}
\newcommand{\sympy}{\textit{SymPy}\xspace}
\newcommand{\latex}{\LaTeX\xspace}
\newcommand{\autofreefem}{\textit{AutoFreeFem}\xspace}

\newcommand{\ie}{\textit{i.e.}\,}
\newcommand{\eg}{\textit{e.g.}\,}
\newcommand{\R}{\mathbb R}

\newcommand{\dd}{\;\mathrm{d} }
\newcommand{\ddx}{\;\mathrm{d} {\mathbf x}}
\newcommand{\ddsx}{\;\mathrm{d} s_{\mathbf x}}
\newcommand{\ddX}{\;\mathrm{d} {\mathbf X}}

\newcommand{\pos}{{\mathbf x}}
\newcommand{\posX}{{\mathbf X}}

\newcommand{\objectiveFunction}{J}
\newcommand{\shapeFunction}{\mathcal J}
\newcommand{\shapePerturbation}{T_t}

\newcommand{\euler}[2]{#1^{E,#2}}
\newcommand{\lagrange}[2]{#1^{L,#2}}
\newcommand{\eulerState}{\euler{u}{t}}
\newcommand{\lagrangeState}{\lagrange{u}{t}}
\newcommand{\lagrangeStateAny}{\lagrange{\varphi}{t}}

\newcommand{\OO}[1]{\mathcal{O}(#1)}
\newcommand{\increment}{\eta}
\newcommand{\Frechet}{Fr\'echet\xspace}

\newtheorem{remark}{Remark}%
\newtheorem{model}{Model}%

\definecolor{codegreen}{rgb}{0,0.6,0}
\definecolor{codegray}{rgb}{0.5,0.5,0.5}
\definecolor{codepurple}{rgb}{0.8,0,0.5}
\definecolor{backcolour}{rgb}{0.9,0.9,0.9}

\lstdefinestyle{mystyle}{
	backgroundcolor=\color{backcolour},   
	commentstyle=\color{codegreen},
	keywordstyle=\color{magenta},
	numberstyle=\tiny\color{codegray},
	stringstyle=\color{codepurple},
	basicstyle=\ttfamily\footnotesize,
	breakatwhitespace=false,         
	breaklines=true,                 
	captionpos=b,                    
	keepspaces=true,                 
	numbers=left,                    
	numbersep=5pt,                  
	showspaces=false,                
	showstringspaces=false,
	showtabs=false,                  
	tabsize=2
}

\pgfplotsset{compat=1.18}

\begin{document}

\title[Article Title]{AutoFreeFem: Automatic code generation with FreeFEM++
	 and LaTex output for shape and topology optimization of non-linear multi-physics problems}


\author[1]{\fnm{Gr\'egoire} \sur{Allaire}}\email{gregoire.allaire@polytechnique.fr}

\author*[1,2]{\fnm{Michael H.} \sur{Gfrerer}}\email{gfrerer@tugraz.at}


\affil[1]{\orgdiv{CMAP, UMR 7641}, \orgname{Ecole Polytechnique, Institut Polytechnique de Paris}, 
	 \city{Palaiseau}, \postcode{91128}, 
	  \country{France}}

\affil[2]{\orgdiv{Institute of Applied Mechanics}, \orgname{Graz University of Technology}, \orgaddress{\street{Technikerstra\ss e 4}, \city{Graz}, \postcode{8010}, \state{Styria}, \country{Austria}}}



\abstract{For an educational purpose we develop the \python package \autofreefem which generates all ingredients for shape optimization with non-linear multi-physics in \freefem and also outputs the expressions for use in \latex. As an input, the objective function and the weak form of the problem have to be specified only once. This ensures consistency between the simulation code and its documentation. In particular, \autofreefem provides the linearization of the state equation, the adjoint problem, the shape derivative, as well as a basic implementation of the level-set based mesh evolution method for shape optimization. 
For the computation of shape derivatives we utilize the	mathematical Lagrangian approach for differentiating PDE-constrained shape functions. Differentiation is done symbolically using \sympy. In numerical experiments we verify the accuracy of the computed derivatives. Finally, we showcase the capabilities of \autofreefem by considering shape optimization of a non-linear diffusion problem, linear and non-linear elasticity problems, a thermo-elasticity problem and a fluid-structure interaction problem. 
}

\keywords{code generation, FreeFEM++, shape derivative, shape optimization}



\maketitle


\section{Introduction}
In order to solve challenging engineering problems, numerical simulation along with shape and topology optimization tools have become an integral part of the design process. Since the computation of linearizations and error-prone shape derivatives for non-linear multi-physics problems are involved, we have developed an educational tool for their automatic code generation. The developed tool outputs always two representations of each expression: one representation for producing a \latex documentation and one representation for writing a simulation and shape optimization script in \freefem.

\subsection{Principles and used software}
\autofreefem is an open-source \python package and can be  downloaded at \href{https://gitlab.tugraz.at/autofreefem/autofreefem}{https://gitlab.tugraz.at/autofreefem/autofreefem}. All example files, discussed in this paper, can be found on this repository. 
It builds on the symbolic calculation capabilities of the open-source \python library \sympy \cite{sympy} (see \Cref{sec::sympy} for details).
For the numerical simulation and shape/topology optimization, the popular open-source software \freefem \cite{hecht2012new} is utilized. \freefem is designed for the efficient numerical solution of partial differential equations using the finite element method in both two and three dimensions. For the  documentation of the problem (input equations, linearization, adjoint problem, shape derivative), the typesetting system \latex is used.  
Thus, the main philosophy of \autofreefem is to provide an implementation in a \freefem script and a documentation in \latex from a single source and therefore allows the fast and reliable development of solutions to complicated problems. In order to illustrate this principle, we consider the following elementary example.
\paragraph{Example: divergence of a vector field}
	Consider the divergence of a vector field which is denoted by the symbol $u$. \Cref{ex::divergenceVectorTab} gives the corresponding outputs  for  \latex and \freefem, respectively. In order to distinguish a vector field from a scalar field, it is printed in bold font in \latex. On the other hand in \freefem, we need to define components, \ie $u_x,\,u_y$ for a 2d problem and $u_x,\,u_y,\,u_z$ in 3d. 
	\begin{table*}[ht]
		\begin{center}
			\begin{tabular}{ccc}
				\toprule
				\latex & \freefem (2D) & \freefem (3D) \\
				\textbackslash operatorname\{div\} \textbackslash mathbf\{u\} & dx(ux) + dy(uy) & dx(ux) + dy(uy) + dz(uz)\\
				\bottomrule
			\end{tabular}
		\end{center}
		\caption{\latex and \freefem expressions for the divergence of a vector field}
		\label{ex::divergenceVectorTab}
	\end{table*}
	\autofreefem works with a versatile unified input that dynamically adapts to various use cases. In the current example, the input takes the form of 
 \begin{equation*}
     \text{div(VectorField('u', ...)}
 \end{equation*}
 Here, the class \textbf{div} (see \Cref{sec::DifferentialOperators}) implements the divergence, whereas the class \textbf{VectorField} is used to define the vector field $\mathbf u$. It is worth noting that the vector field, although not fully depicted here for simplicity, requires five input arguments, which give information on the domain of definition and the boundary conditions (see \Cref{sec::fields} for details). 
 
\paragraph{Automatic simulation and shape optimization in \freefem}\label{ex::simulationOptimization} 
The main class of \autofreefem is called \textbf{Lagrangian} (see \Cref{sec::fields}). This class has in particular the two methods 'setup\_simulation' and 'setup\_optimization'. When the first method is called, \autofreefem first checks if the problem is linear or non-linear (see \Cref{sec::Classsification}). In case of a linear problem, the output is a simple \freefem script for the simulation of the problem. In case of a non-linear problem, the linearization is computed and used in a Newton's method implemented in the output \freefem script. For simulations the input of the Lagrangian are the primary field variables, the corresponding test functions and a variational formulation of the problem (see \Cref{sec::flow}). For shape optimization problems the method 'setup\_optimization' has to be called. Then a linearization in case of a non-linear problem, an adjoint problem, and a shape derivative are computed. Furthermore, a corresponding \freefem script is generated (see e.g. \Cref{sec::linearCantilever}). For the numerical solution of shape optimization problems in \freefem we employ the level-set based mesh evolution method \cite{allaire2014shape}. To this end we use the aditional open-source libraries \textit{mmg}\footnote{\href{http://www.mmgtools.org/}{http://www.mmgtools.org/}} \cite{dapogny2014three}, \textit{mshdist}\footnote{\href{https://github.com/ISCDtoolbox/Mshdist}{https://github.com/ISCDtoolbox/Mshdist}} \cite{dapogny2012computation}, and \textit{advection}\footnote{\href{https://github.com/ISCDtoolbox/Advection}{https://github.com/ISCDtoolbox/Advection}}.

\subsection{Relation to the literature}
The automatic generation of simulation code and the automatic computation of shape derivatives have been considered in some previous works. As part of the FEniCS Project \cite{alnaes2015fenics}, the Unified Form Language (UFL) offers a flexible interface for choosing finite element spaces and defining expressions for weak forms in a notation close to mathematical notation \cite{alnaes2014unified}. This allows also for the automatic computation of derivatives and therefore the easy treatment of non-linear problems. Based on UFL, the open-source library FEMorph is an automatic shape differentiation toolbox, which can compute first- and second-order shape derivatives \cite{schmidt2018weak}. It refactors UFL expressions and applies shape calculus differentiation rules recursively. In \cite{ham2019automated}, the UFL is extended to shape differentiation using a different strategy. The approach in \cite{ham2019automated} is based on pullbacks and standard Gateaux derivatives. Furthermore,  automated shape derivatives for transient PDEs in FEniCS and Firedrake \cite{rathgeber2016firedrake} are presented in \cite{dokken2020automatic}. This has been further developed in the software Fireshape \cite{paganini2021fireshape}. Inspired by the FEniCS Project, the finite element software package NGSolve \cite{schoberl2014c} has a flexible interface to \python, which allows  defining expressions for weak forms in a mathematical notation. In \cite{gangl2021fully} it has been extended for the automatic computation of first- and second-order shape derivatives based on a Lagrangian function, pullbacks, and directional (Gateaux) derivatives. In \cite{gangl2022automated}, NGSolve has been further developed to allow also for the automatic computation of topological derivatives. We mention also the software \textit{cashocs} described in \cite{BLAUTH2021100646,blauth2023version}, which offers automated solutions for shape optimization and optimal control. 

All works mentioned so far are based on symbolic shape differentiation. In the context of density based methods for topology optimization, an automatic differentiation (AD) tool is presented in \cite{chandrasekhar2021auto}. There, practically no difference in the timings of AD and symbolic sensitivities was found. For further references on AD we refer to \cite{chandrasekhar2021auto}.
In all these references the focus is on the automatic generation of the sensitivity information for use within the computational optimization routine. 
The main novelty of the present work is to consider, in a pedagogical perspective, the simultaneous generation of \latex expressions for the documentation and \freefem expressions for numerical optimization. 

\subsection{Outline of the paper}
The next section is a brief tutorial which features two examples: the simulation of a 3d non-linear fluid flow and the 2d shape optimization of a linearly elastic structure. 
In \Cref{sec::Theory} we present the underlying mathematical theory of \autofreefem. \Cref{sec::implementation} details the implementation. Several numerical examples are discussed in \Cref{sec::examples}. Finally, we draw some conclusions from the present work in \Cref{sec::conclusion}.

\section{Introductory examples}
This section provides two hands-on introductory examples to \autofreefem. First, the simulation of a viscous fluid flowing through a pipe. Second, the compliance minimization of a cantilever beam.

\subsection{3d simulation of a fluid flow}\label{sec::flow}
In this example, the fluid flow in a winding pipe is simulated by solving the incompressible Navier-Stokes equations with Taylor-Hood finite elements. In the next subsections, we explain step by step how to solve this problem using \autofreefem.

\subsubsection{Step 1: 3d mesh generation with \freefem}
For the generation of the computational mesh we use built-in commands of \freefem in the file "meshNS3d.txt" (which can be found on the \autofreefem repository). In particular, we use \textit{border} and \textit{buildmesh} to generate a disk, which is then extruded to a cylinder by \textit{buildlayers}. The final mesh is obtained by a mesh distortion using the command \textit{movemesh}. The chosen 3d geometry of a winding pipe is depicted in \Cref{fig::meshFlowA}.
Note that more complicated 3d meshes can be created with another mesh generator and loaded into \freefem. 

\subsubsection{Step 2: Definition of the problem}
The input for \autofreefem for this example is given in \Cref{list::flow}.	
The first line imports all modules from the package \autofreefem. In the lines 4-6 the physical constants $\rho$ (density) and $\mu$ (viscosity), as well as a penalty parameter $\gamma$ are defined using the class \textbf{Constant}. Next, in line 10, we use the class \textbf{VectorField} to define the fluid velocity $\mathbf v$, which is discretized by finite elements of polynomial degree 2 (P2) on the domain \textit{Th}. Furthermore, here Dirichlet boundary conditions (prescribed velocities) on boundaries with labels 1, 3 and 4 are also defined. In order to define the corresponding vector-valued inhomogeneous boundary data, we use the class \textbf{BoundaryFunction} in line 9. In line 11 we define a corresponding test function $\delta \mathbf v$ to the velocity field $\mathbf v$. In lines 13 and 14 the class \textbf{ScalarField} is used to define the fluid pressure $p$ discretized by P1 finite elements and the test function $\delta p$. A homogeneous pressure is prescribed at the boundary with label 2. In order to increase the readability in the \latex output, we introduce the viscous fluid stress $\pmb \sigma_{f}$ as an \textbf{Expression} in lines 16 and 17 (see also \Cref{sec::fields}). In lines 19-25 a \textbf{Lagrangian} object is set up. The first argument is a list of all primary field functions $(\mathbf v, p)$, whereas the second argument $(\delta \mathbf v, \delta p)$ is a list of all test functions. The third argument is the objective function, which has no meaning in case of a simulation. The fourth argument is the weak formulation of the governing equations. Those are the momentum balance and the continuity equation enhanced by a penalty term. For example, the input in line 20 
represents the viscous stress term in the momentum equation and the corresponding \latex representation is
\begin{equation}
    \int_{\Omega} \left( \pmb \sigma_{f}( {\mathbf{ v } }) \mathrel{:}  \nabla \delta{\mathbf{ v}}  \right) \,dx.
\end{equation}
Note, that the domain integral over $\Omega$ is realized by the class \textbf{dx} (see \Cref{sec::Integrals} ; not to be confused with the $x$-derivative notation in \freefem), whereas \textbf{inner2} (see \Cref{sec::TensorAlgebra}) implements the double dot product.

Finally, in line 26, we call the method \textit{setup\_simulation} of the Lagrangian and specify the file "meshNS3d.txt" for the definition of the finite element mesh.
\begin{lstlisting}[language=Python,caption={Input file for the Navier-Stokes example (\href{https://gitlab.tugraz.at/autofreefem/autofreefem/-/blob/main/examples/simulation/run_navierStokes.py}{run\_navierStokes.py})},captionpos=b,label={list::flow}]
from autofreefem import *

name = 'Navier-Stokes'
rho = Constant('\\rho', 1000., 'fluid density', 'kg/m^3')
mu = Constant('\mu', 1.0e-3, 'fluid viscosity', 'N s/m^2')
penalty = Constant('\\gamma', 1.e-9, 'penalty term', '1/(Pa s)')
Th = Domain('\\Omega','Th')
# fluid velocity
diri = BoundaryFunction( '0', '0', '5e-2*(r^2-x^2-y^2)*(z<0.05)')
v = VectorField('v', 'P2', Th, diri, '1/3/4')
testv = VectorField('\delta v', 'P2', Th, '0.', 'none')
# fluid pressure
p = ScalarField('p', 'P1', Th, '0.', '2')
q = ScalarField('\delta p', 'P1', Th, '0.', 'none')
# viscous fluid stress tensor
d = (grad(v) + transpose(grad(v))) / 2
sigmaFluidViscous = Expression('\pmb \sigma_f', mu * d)
#
lag = Lagrangian([v, p], [testv, q], 0,
    dx((inner2(sigmaFluidViscous, grad(testv))), Th)  # viscous stress
    - dx(div(testv) * p, Th)   # pressure
    + dx(rho * inner(inner(grad(v), v), testv), Th)  # convection
    + dx(penalty * p * q, Th)  # penalty term
    - dx(div(v) * q, Th), # continuity equation
    dimensions=3)
lag.setup_simulation(name, mesh='meshNS3d.txt')


\end{lstlisting}
\subsubsection{Step 3: \latex output}
The simulation problem defined in \Cref{list::flow} is documented using \latex. The output is given in \Cref{fig::FlowDocu}. The output starts with the user input of the governing equations (here the momentum equation and the continuity equation). Then the abbreviations used, \ie all objects of class \textbf{Expression}, are defined (here the viscous fluid stress $\pmb \sigma_{f}$). \autofreefem automatically detects that the problem is non-linear (due to the here non-linear convective term, see \Cref{sec::Classsification}) and computes also a linearization for the use within a Newton method. Thus, the problem for the Newton update is given in the remainder of \Cref{fig::FlowDocu}. In addition, the numerical values of the considered physical parameters are supplied in the automatically generated \Cref{tab::FlowParameter}.
\begin{figure*}[htb]
	\fbox{\parbox{\textwidth}{
			\input{NavierStokes.tex}
	}}
	\caption{Automatically generated \latex documentation of the problem formulation and linearization of the Navier-Stokes problem defined in \Cref{list::flow}.}
	\label{fig::FlowDocu}
\end{figure*}
\begin{table}[htb]
	
\begin{spacing}{1.2}
\begin{tabular}{cccl}
\toprule
fluid density & $\rho$ & 1.00e+03 & $kg/m^{3}$ \\ 
fluid viscosity & $\mu$ & 1.00e-03 & $N s/(m^{2)}$ \\ 
penalty term & $\gamma$ & 1.00e-09 & $1/(Pa  s)$ \\ 
\bottomrule
\end{tabular}
\end{spacing}

	\caption{Automatically generated \latex documentation of numerical values of constants defined in \Cref{list::flow}}
	\label{tab::FlowParameter}
\end{table}
\subsubsection{Step 4: Simulation with \freefem}
In addition to the \latex output, \autofreefem, run in the simulation mode (method 'setup\_simulation'), produces three ".edp" files:
\begin{itemize}
	\item run\_NavierStokes.edp
	\item NavierStokesResidual.edp
	\item NavierStokesNewton.edp
\end{itemize}
The file "run\_NavierStokes.edp" implements a basic solver based on Newton's method for the simulation of the flow problem in \freefem. To this end, it uses the expressions (varf's) defined in the other two files for evaluating the residual vector and the Jacobian matrix. The latter files can also be used as building blocks for more advanced solvers (\eg preconditioned iterative solvers and/or domain decomposition methods) implemented by the user.
\subsubsection{Step 5: Results}
\begin{figure*}
	\begin{subfigure}[t]{0.28\textwidth}
		\cropped{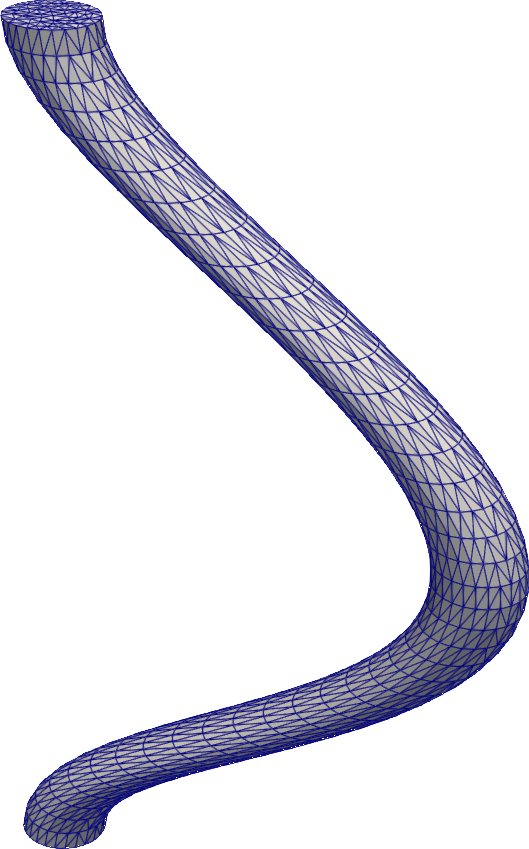}
		\caption{}
		\label{fig::meshFlowA}
	\end{subfigure}
	\hfill
	\begin{subfigure}[t]{0.28\textwidth}
		\cropped{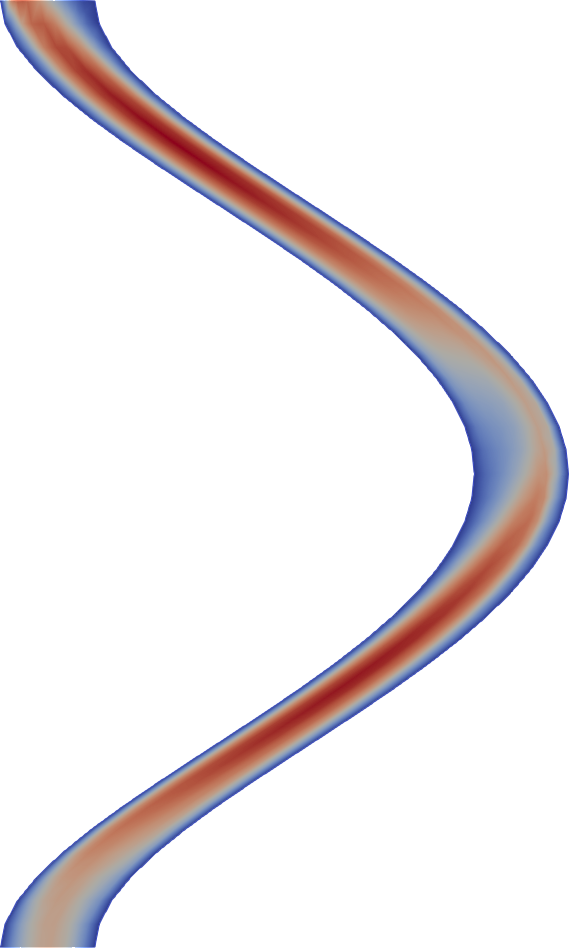}
		\caption{}
		\label{fig::meshFlowB}
	\end{subfigure}
	\hfill
	\begin{subfigure}[t]{0.28\textwidth}
		\cropped{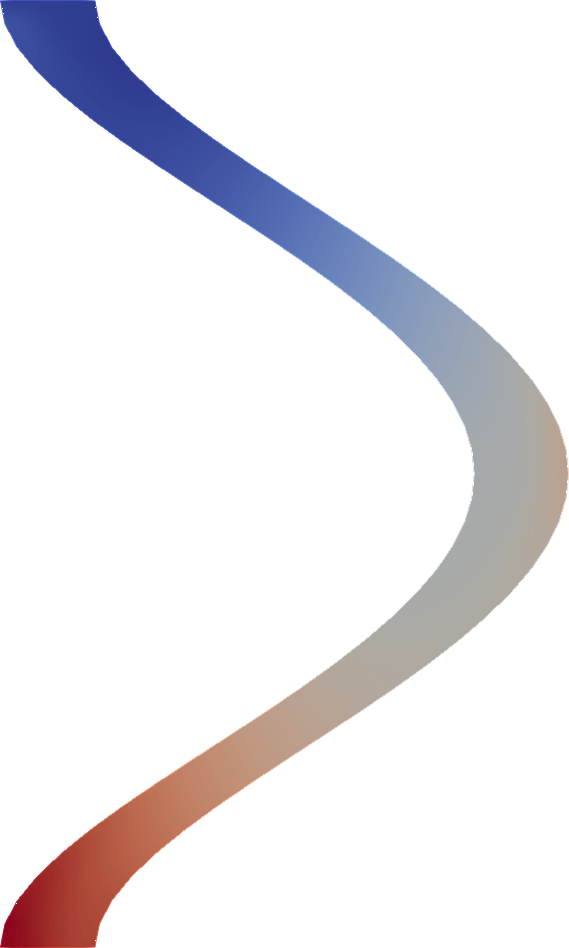}
		\caption{}
		\label{fig::meshFlowC}
	\end{subfigure}
	\caption{Simulation of the fluid flow in a
pipe: (a) 3d computational mesh; (b) Distribution of the norm of the fluid velocity over a vertical slice of the domain. Red
		corresponds to high velocity, blue corresponds to low velocity. On the lower surface the velocity distribution is prescribed as Dirichlet boundary condition. (c) Distribution of the (relative) fluid pressure over a vertical slice of the domain. Red
		corresponds to high pressure, blue corresponds to low pressure. The pressure on the upper surface is prescribed as Dirichlet boundary condition.  }
	\label{fig::meshFlow}
\end{figure*}

Running the file "run\_NavierStokes.edp" with \freefem produces a ".vtu" file with the simulation results. This file can be viewed with a 3D graphics program such as \textit{paraview}\footnote{\href{https://www.paraview.org/}{https://www.paraview.org/}}. The computed fluid velocity and pressure are depicted in \Cref{fig::meshFlow}.
\subsection{Shape optimization of a cantilever}
\label{sec::linearCantilever}
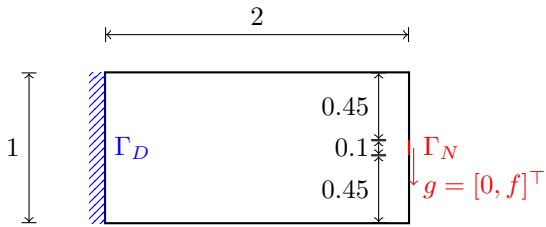
\begin{figure}[ht]
	\centering
	\begin{tikzpicture}[scale=2.]
		\useasboundingbox (-0.7,-0.3) rectangle (3.4,1.5);
		\draw[thick] (0,0) rectangle (2,1);
		\fill[pattern=north east lines,pattern color=blue] (-0.1,0) rectangle (0,1);
		\draw[blue] (0,0) -- (0,1) node[pos=0.5,right]{$\Gamma_D$};
		\draw[->,red] (2.03,0.5) node[right]{$\Gamma_N$} --++(0,-0.25) node[right]{$g = [0,f]^\top$} ;
		\draw[red,thick] (2,0.45) --++(0,0.1);
		\draw[|<->|] (0,1.25) --++(2,0) node[pos=0.5,above]{2};
		\draw[|<->|] (-0.5,0.) --++(0.0,1) node[pos=0.5,left]{1};
		\draw[|<->|] (1.8,0.) --++(0.0,0.45) node[pos=0.5,left]{0.45};
		\draw[|<->|] (1.8,0.45) --++(0.0,0.1) node[pos=0.5,left]{0.1};
		\draw[|<->|] (1.8,0.55) --++(0.0,0.45) node[pos=0.5,left]{0.45};
	\end{tikzpicture}
	\caption{Geometry and boundary conditions of the elastic cantilever. }
	\label{fig::Cantilever}
\end{figure}

Our second introductory example is a classical 2d compliance minimization for an elastic cantilever.
The working domain is a rectangle of size $2 \times 1$ with zero displacement boundary condition on the left side and a vertical load applied on a small portion of length $0.1$ at the middle of the right side, denoted by $\Gamma_N$, such that the resultant force has unit magnitude. All other boundaries are traction free. The geometry and the boundary conditions are illustrated in \Cref{fig::Cantilever}. There are no body
forces.
\subsubsection{Step 1: Mesh}
Again, for the generation of the computational mesh, we use built-in commands of \freefem in the file \href{https://gitlab.tugraz.at/autofreefem/autofreefem/-/blob/main/examples/shapeOptimization/elasticity/meshelasticCantilever.txt}{meshelasticCantilever.txt}.

\subsubsection{Step 2: Definition of the problem}
\begin{lstlisting}[language=Python,caption={Input file for the cantilever optimization (\href{https://gitlab.tugraz.at/autofreefem/autofreefem/-/blob/main/examples/shapeOptimization/elasticity/run_linearElasticity.py}{run\_linearElasticity.py})},captionpos=b,label={list::Cantilever}]
from autofreefem import *
name = 'linear Elasticity'
bndN = Boundary('\\Gamma_N','ThL','2')
Th = Domain('D','Th',bndN)
# definition of the displacement field and the test function
u = VectorField('u', 'P2', Th, '0.', '4')  # left fixed
testu = VectorField('\delta u', 'P2', Th, '0.', 'none')
# definition of constants
f = Constant('f', -10, 'vertical load component', 'N/m')
nu = Constant('\\nu', .3, 'Poisson\'s ratio', '-')
E = Constant('E', 200., 'Young\'s modulus', 'N/m^2')
Afac = Constant('\ell', 0.25, 'Lagrange multiplier', '-' )
# definition of the material contrast
X = CharacteristicFunction('\\raisebox{\\depth}{$\\chi$}', 1, '1/100')
# definition of expressions
mu = Expression(sp.Symbol('\mu'), E/(2*(1+nu)))
lam = Expression(sp.Symbol('\lambda'), E*nu/((1-2*nu)*(1+nu)))
strain = Expression(sp.Symbol('\pmb \\varepsilon'), (grad(u) + transpose(grad(u))) / 2)
stress = Expression(sp.Symbol('\pmb \\sigma'), X*(lam*tr(strain)*identity() + 2 * mu*strain) )
# definition of the shape optimization problem
J = Expression('J', dx(inner2(stress,strain), Th) + dx(Afac*X, Th))
lag = Lagrangian([u], [testu], J, dx(inner2(stress, grad(testu)), Th) - dsx(f*inner(testu, ey()), Th, 0))
lag.setup_optimization(name, mesh='meshelasticCantilever.txt',             
        hmin=0.01, hmax=0.02, diffusion=1./10000., v0=0.002, iterations=82,
        phi='-0.4 - sin(pi * kx * (x+0.5)) * cos(pi * ky * (y))',
        boundaryLabels=['1','4'], fixedLabels=['2'], show_weak_material=False)
\end{lstlisting}
The \python input for \autofreefem is given in \Cref{list::Cantilever}. In lines 1-12 the displacement field, the physical constants and the used expressions are defined. These commands were already used in \Cref{sec::flow}. Thus, we focus on the following lines which involve new aspects due to the considered shape optimization problem. In line 14 we introduce an object of the class \textbf{CharacteristicFunction} in order to distinguish between solid and void material. The first input is a \latex symbol ({\raisebox{\depth}{$\chi$}}) for this function. The second argument represents the relative strength of the solid material (typically 1), whereas the third argument refers to the void material. Here, we use a factor of $\nicefrac{1}{100}$ in order to mimic void by a very weak material. Strictly speaking, the object \textbf{CharacteristicFunction}, defined in line 14, is not the characteristic function of a set but rather a "color function", taking two different values (not necessarily 0 and 1) in two sub-domains. In lines 16 to 19 we use the class \textbf{Expression} to define the Lam{\'e} constants, the strain tensor and the stress tensor. In line 21 the objective function $J$ is defined as a combination of the compliance and a fixed Lagrange multiplier $\ell$ multiplied by the area of the solid. Next, in line 22 the \textbf{Lagrangian} object is set up. The load on the boundary is incorporated in the problem by the corresponding boundary integral using the class \textbf{dsx} (see \Cref{sec::Integrals}). Finally, in lines 22-26 we call the method \textit{setup\_optimization} of the Lagrangian and specify the file "meshelasticCantilever.txt" for the definition of the finite element mesh. Specifically for shape optimization, we also specify a minimal (\textit{hmin}) and maximal (\textit{hmax}) element size for remeshing with \textit{mmg} \cite{dapogny2014three}, a diffusion parameter for the regularization in the shape gradient identification problem \cite{allaire2021shape}, an initial optimization velocity $v0$ used in the advection of the level-set function, a maximal iteration number and an initial level-set function $\phi$ defining the initial design. Furthermore, we set labels for boundaries where the normal component of the shape gradient should be set to zero (\textit{boundaryLabels}), and labels for boundaries where the full gradient should vanish (\textit{fixedLabels}). Finally, we set the postprocessing option \textit{show_weak_material} such that the weak material is not shown in the ".vtu" outputs.

\subsubsection{Step 3: \latex output}
The problem formulation defined in \Cref{list::Cantilever} is documented in Appendix \ref{appendix::Cantilever}. \autofreefem automatically detects that the state problem is linear and therefore skips the statement of a superfluous linearization. Furthermore, the numerical values of the considered physical parameters are supplied in the automatically generated \Cref{tab::CantileverParameter}.
\begin{table}[htb]
	
\begin{spacing}{1.2}
\begin{tabular}{cccl}
\toprule
Lagrange multiplier & $\ell$ & 0.25 & $\frac{N}{m^2}$ \\ 
Poisson's ratio & $\nu$ & 0.3 & $-$ \\ 
Young's modulus & $E$ & 200 & $\frac{N}{m^2}$ \\ 
vertical load component & $f$ & -10 & $\frac{N}{m^2}$ \\ 
\bottomrule
\end{tabular}
\end{spacing}
	\caption{Automatically generated \latex documentation of numerical values of constants defined in \Cref{list::Cantilever}}
	\label{tab::CantileverParameter}
\end{table}

\subsubsection{Step 4: Simulation with \freefem}
In addition to the \latex output, \autofreefem in the shape optimization mode (method 'setup\_optimization') produces six ".edp" files:
\begin{itemize}
	\item optimize\_linearElasticity.edp
	\item linearElasticityObjective.edp
	\item linearElasticityResidual.edp
	\item linearElasticityNewton.edp
	\item linearElasticityAdjoint.edp
	\item linearElasticityShape.edp
\end{itemize}
The file "optimize\_linearElasticity.edp" implements a basic solver for the simulation of the cantilever and a basic variant of the level-set based mesh evolution method introduced in  \cite{allaire2014shape,allaire2021shape}. To this end, it uses the expressions (varf's) defined in the other files for evaluating the objective functional, the residual vector, the stiffness matrix, the resolution of the adjoint problem and finally the shape derivative.
\subsubsection{Step 5: Results}
\begin{figure}[ht]
	\begin{subfigure}{0.49\textwidth}
		\cropped{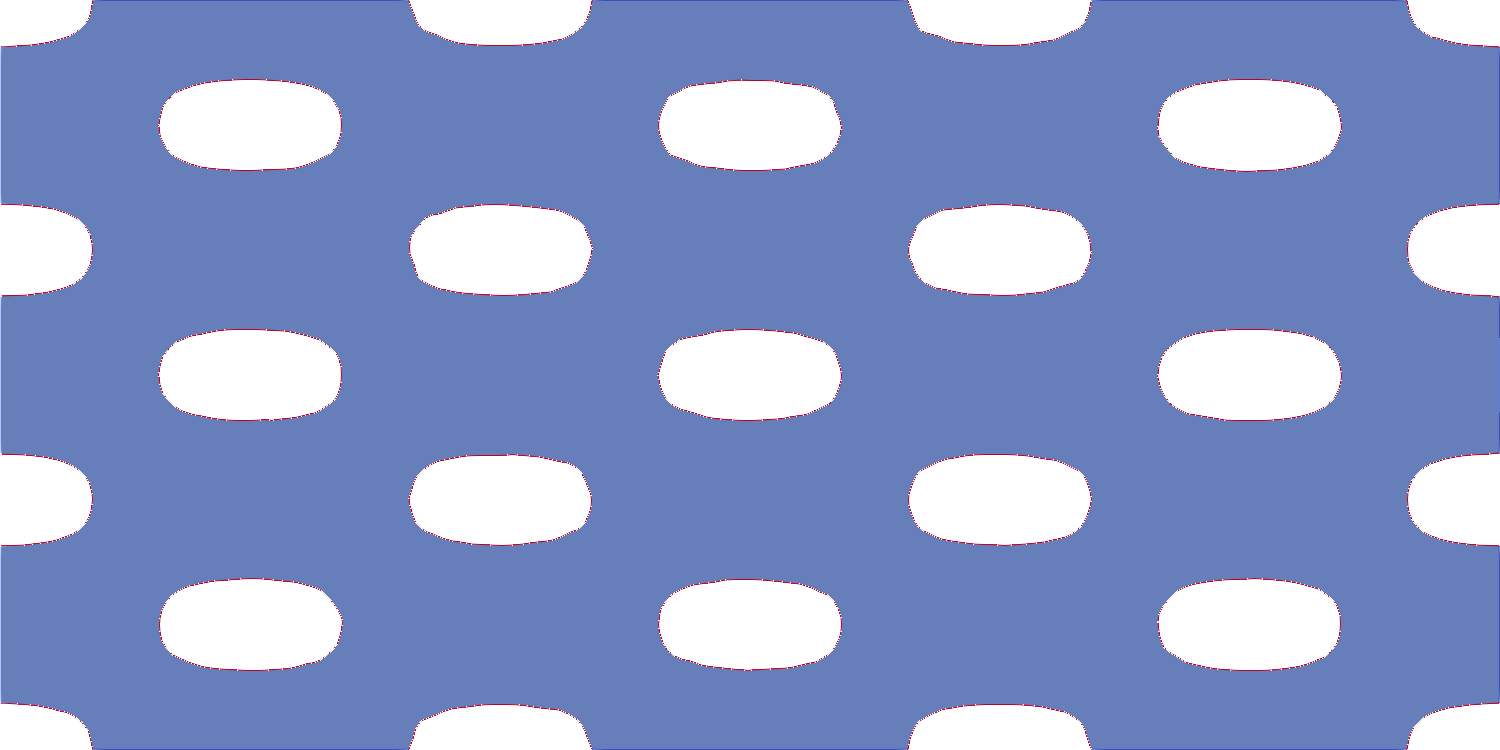}
		\caption{}
		\label{fig:picCantileverElasticInitialization}
	\end{subfigure}
	\begin{subfigure}{0.49\textwidth}
		\cropped{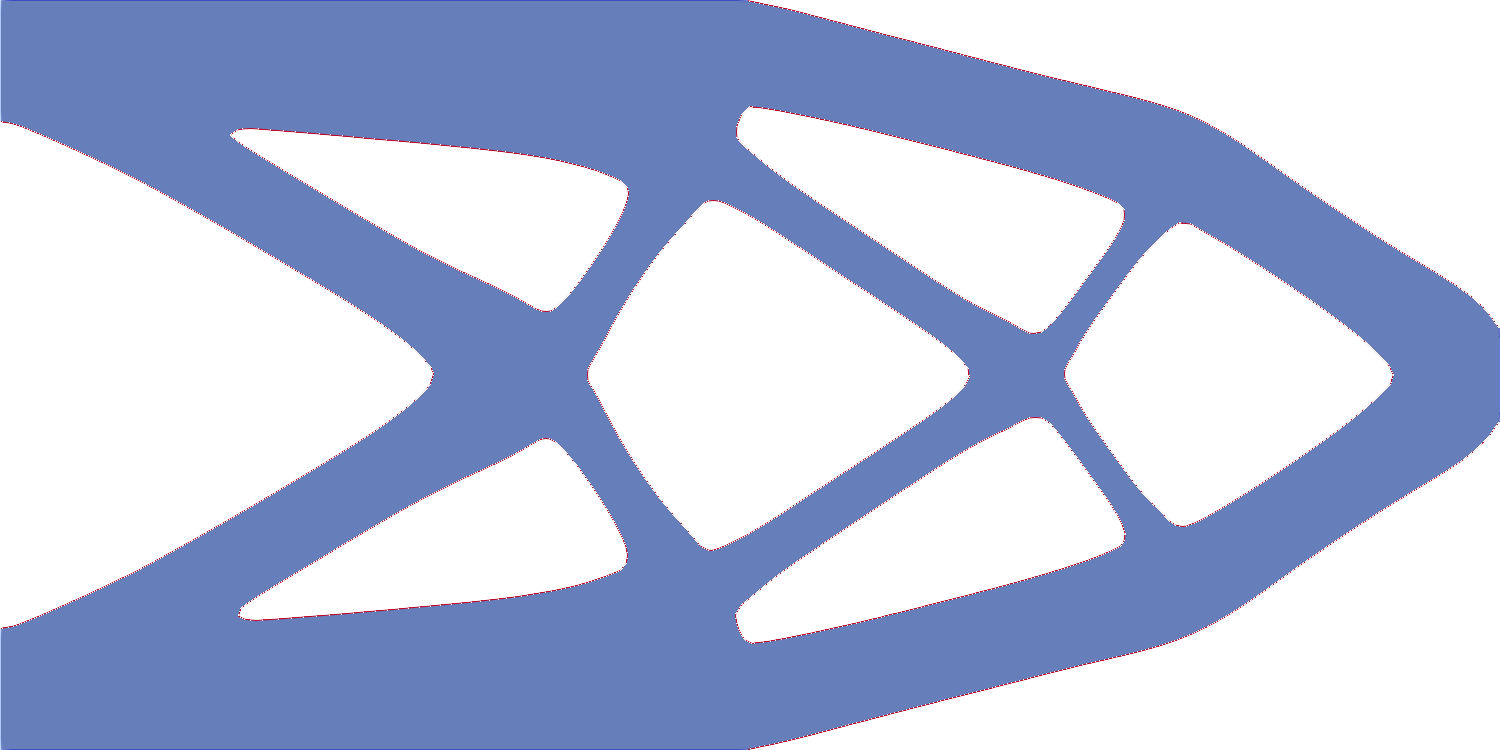}
		\caption{}
		\label{fig:picCantileverElasticOptimal}
	\end{subfigure}
	\caption{Linear elastic cantilever problem: (a) Initialization; (b) Optimized design. }
	\label{fig:picCantileverElastic}
\end{figure}

\begin{figure}[ht]
	\centering
	\begin{tikzpicture}
		\begin{axis}[width=0.8\linewidth,
			grid=both,
			grid style={line width=.1pt, draw=gray!10},
			major grid style={line width=.2pt,draw=gray!50},
			minor tick num=5,
			legend cell align=left,
			legend pos=outer north east,
			legend style={draw=none},
			ylabel=$J$,
			xlabel=iteration]
			\addplot +[mark minAbove, every node near coord/.style=] table [x=iteration,y=objective,col sep=comma] {objectivelinearElasticity.txt};
		\end{axis}
	\end{tikzpicture}
	\caption{Convergence history for the linear elastic cantilever problem}
	\label{fig:convergenceElasticCantilever}
\end{figure}
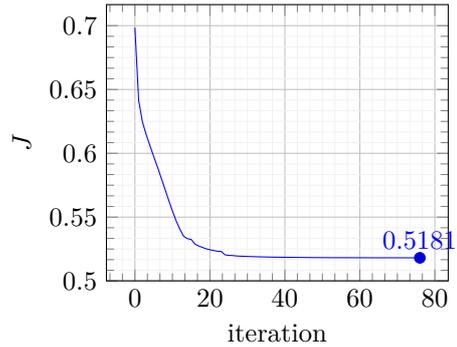

Running the file "optimize\_linearElasticity.edp" with \freefem produces a sequence of 200 ".vtu"-files with the optimization results. The initialization and the optimized design are depicted in \Cref{fig:picCantileverElastic}.  The evolution of the objective function is reported in \Cref{fig:convergenceElasticCantilever}.
\section{Theory on optimization problems involving non-linear multi-physics PDE constraints}\label{sec::Theory}
In this section, we briefly describe the underlying mathematical theory of \autofreefem. In the present paper, we consider PDE-constrained shape optimization problems of the form:
\begin{align}\label{eq::abstractOptimization}
	\begin{aligned}
		\min_{\Omega \in \mathcal A}\;& \objectiveFunction(\Omega,u)\\
		&\text{subject to  $u\in V(\Omega)$:}\\
		&R(\Omega,u,v) = 0 \quad \mbox{for all } v \in V(\Omega). 
	\end{aligned}	
\end{align}
Here, $\mathcal A$ is the  a set of admissible shapes, $\objectiveFunction:\mathcal A \times V(\Omega) \rightarrow \R$ the objective function and the state $u$ is the solution of the governing non-linear physics incorporated in $R(\Omega,u,v)$ posed on a Hilbert space $V(\Omega)$. For given shape $\Omega\in \mathcal A$ let $u(\Omega)$ be the unique solution of the state equation. This allows to introduce the shape functional $\shapeFunction:\mathcal A \rightarrow \R$,  
\begin{equation}
	\shapeFunction(\Omega) =  \objectiveFunction(\Omega,u(\Omega)).
\end{equation}
In \Cref{sec::Classsification} and \Cref{sec::Linearization}, we first focus on the classification and resolution of the possibly non-linear state equation $R(\Omega,u,v) = 0$. Following this, we describe the theory on computing the shape derivative for single-physics problems in \Cref{sec::Single} and for multi-physics problems in \Cref{sec::Multi}.

\subsection{Classifying the State Problem: Linear vs. Non-Linear}\label{sec::Classsification}
In a linear state problem, $R(\Omega,u,v)$ contains only terms that are either independent of the solution field $u$ or linearly dependent on $u$. So  we check whether the equation 
\begin{equation}\label{eq::classification}
    \frac{d^2}{d\tau^2}R(\Omega, \tau u, v) = 0
\end{equation}
is satisfied or not. Indeed, \eqref{eq::classification} holds for a linear state problem, but not for a nonlinear problem. In order to implement the automatic evaluation of \eqref{eq::classification} and subsequent expressions in this section, we use the \sympy commands \textit{subs} and \textit{diff} (see \Cref{sec::sympy}).
For multi-physics problems (see \Cref{sec::Multi}) condition \eqref{eq::classification} generalizes to 
\begin{equation*}
\frac{d^2}{d\tau^2} R_i(\Omega,\tau u_1,\mydots, \tau u_N, v_i)=0, \text{ for } i=1,\mydots,N.
\end{equation*}
\subsection{Linearization}\label{sec::Linearization}
The solution of a non-linear state problem is typically computed by Newton's method. Thus, at a discrete level, we should provide the Jacobian matrix after the discretization of $R(\Omega,u,v)$. However, here we compute the expression of the Jacobian matrix at the continuous level. Thus, we follow a \textit{first-differentiate-then-discretize} approach in the present paper, which might not give the same Jacobian matrix obtained by the \textit{first-discretize-then-differentiate} approach. However, for linearizations both approaches typically give the same Jacobian matrix (a typical counterexample is that of plasticity problems, see \cite{wriggers2008nonlinear}). Note, that for the shape derivative the result of \textit{first-differentiate-then-discretize} is usually different from that of \textit{first-discretize-then-differentiate}.

In order to illustrate the abstract setting in \eqref{eq::abstractOptimization}, consider a non-linear diffusion problem with homogeneous boundary conditions on a domain $\Omega$. 
\begin{model}\label{problem::model}
	Find the solution $u\in H^1_0(\Omega)$ such that $R(\Omega,u,v)=0$ for all test functions $v\in H^1_0(\Omega)$, where
    \begin{align}
		R(\Omega,u,v) &= a(\Omega,u,v) - b(\Omega,v), \\ 
    \nonumber a(\Omega,u,v) &= \int_\Omega \lambda(\posX,u(\posX)) (\nabla u \cdot \nabla v) \ddX, \\
    \nonumber b(\Omega,v) &= \int_\Omega f \,v \ddX,
    \end{align}
	with a spatially varying diffusion coefficient $\lambda$, which might depend on the state $u$.
\end{model}
As usual $R(\Omega,u,v)$ is linear with respect to the test function $v$. 
In order to linearize $u\mapsto R(\Omega,u,v)$, we introduce its \Frechet derivative as the linear and bounded form $d_FR(\Omega,u,v)$, which satisfies
\begin{widetext}
\begin{equation}
\begin{aligned}
	R(\Omega,u+\increment,v) = R(\Omega,u,v) + d_FR(\Omega,u,v) \cdot \increment + o(\increment), \quad\text{with} \quad \lim_{\increment \rightarrow 0} \frac{|o(\increment)|}{\|\increment\|} = 0.
\end{aligned}    
\end{equation}    
\end{widetext}
Then, the linearization $LR$ of $R$ at state $u_0$ reads \cite{wriggers2008nonlinear}
\begin{align*}
	LR(\Omega,u_0,v,\increment) &= R(\Omega,u_0,v) + d_FR(\Omega,u_0,v)\cdot\increment ,
\end{align*}
where $LR$ is linear in the third and the last argument, but still non-linear in the first two arguments. In Newton's method, the update of a known state $u_0$ reads $u_0 \leftarrow u_0 + \increment$, where $\increment$ is the solution to $LR(\Omega,u_0,v,\increment) = 0$, which has to hold for all test functions $v\in V$. In order to practically compute the linearization, it is more convenient to use the notion of Gateaux (or directional) derivative 
\begin{equation}
	\begin{aligned}\label{eq::Gateaux}
		dR(\Omega,u,v,\increment) &=  \frac{d}{d\tau}R(\Omega,u+\tau\increment,v)\bigg|_{\tau=0}.
	\end{aligned}	
\end{equation} 
If $R$ is \Frechet differentiable, then it is Gateaux differentiable too and $d_FR=dR$. 
Note that for \Cref{problem::model} we have 
\begin{align*}
	dR(\Omega,u_0,v,\increment) &= \int_\Omega \lambda(\posX, u_0) (\nabla \increment\cdot \nabla v) \ddX \\ \nonumber
	&+ \int_\Omega \frac{\partial\lambda}{\partial u}(\posX,u_0) \increment (\nabla u_0 \cdot \nabla v) \ddX.
\end{align*}
\subsection{Shape derivative of single-physics problems}\label{sec::Single}
In this section, we briefly review the Lagrangian method for computing the shape derivative of a PDE-constraint objective function. Mathematically rigorous treatments can be found in \cite{hinze2008optimization,sturm2015minimax,henrot2018shape}; here we prefer a pedagogical presentation using notations from continuum mechanics. 
The notion of shape derivative relies on Hadamard method, which considers variations of the domain $\Omega\subset \R^d$ of the form
\begin{equation}\label{eq::Hadamard}
	\Omega_t = (Id + t\mathbf V)(\Omega) = T_t(\Omega), 
\end{equation}
where $\mathbf V:\Omega \rightarrow \R^d$ is a vector field and $t$ a scalar perturbation parameter. In the following, we use the notation that $\posX$ is a point in the unperturbed domain $\Omega$ and $\pos$ a point in the perturbed domain $\Omega_t$.
Thus, a point $\posX\in\Omega$ is therefore mapped to $\pos \in \Omega_t$ by
\begin{equation}\label{eq::HadamardPos}
	\pos = T_t(\posX) = \posX + t\mathbf V(\posX),
\end{equation}
and the Jacobian matrix $\mathbf F_t: \Omega \rightarrow \R^{d\times d}$ reads
\begin{equation} \label{eq::Jacobian}
	\mathbf F_t(\posX) = \nabla T_t(\posX) = \mathbf I + t \nabla\mathbf V(\posX),
\end{equation}
where $\mathbf I$ denotes the identity matrix. 
From a continuum mechanics viewpoint, $\posX$ is a Lagrangian coordinate, while $\pos$ is an Eulerian coordinate, 
the vector field $\mathbf V$ can be interpreted as a displacement field and \eqref{eq::Jacobian} defines the associated deformation gradient.
For some shape functional $\shapeFunction(\Omega)$, the shape derivative is defined as
\begin{equation}\label{eq::shapeDerivativeDef}
	DJ(\Omega)(\mathbf V) = \lim_{t\rightarrow0} \frac{\shapeFunction(\Omega_t)-\shapeFunction(\Omega)}{t}. 
\end{equation}
\begin{remark}
	It is also possible to define the shape derivative as the \Frechet derivative of the mapping $\mathbf V \mapsto \mathcal \shapeFunction((\mathbf I + \mathbf V)(\Omega))$ in $\mathbf V = \mathbf 0$ \cite{allaire2021shape}. 
\end{remark}
In the following we introduce the Eulerian and Lagrangian states, which will be employed in the subsequent derivation of an efficient formula for the computation of the shape derivative \eqref{eq::shapeDerivativeDef}.
For a perturbed domain $\Omega_t$ we have the Eulerian state $\eulerState \in V(\Omega_t)$, which satisfies the state equation
\begin{equation}\label{eq::residual}
	R(\Omega_t,\eulerState,v)=0 \quad \forall v \in V(\Omega_t).
\end{equation}
In a next step we introduce the Lagrangian state $\lagrangeState \in V(\Omega)$ defined by the pull-back of the Eulerian state $\eulerState \in V(\Omega_t)$ to the unperturbed domain $\Omega$
\begin{equation}\label{eq::lagrangeEuler}
	\lagrangeState(\posX) = \eulerState \circ T_t(\posX) = \eulerState(\pos).
\end{equation}
Conversely, this allows to write
\begin{equation}\label{eq::EulerLagrange}
	\eulerState(\pos) = \lagrangeState \circ T_t^{-1}(\pos) = \lagrangeState(\posX).
\end{equation}
Obviously, for $t=0$ in \eqref{eq::HadamardPos}, we have $\pos = \posX$ and thus the Lagrangian and the Eulerian states coincide $u^{L,0} = u^{E,0} = u^0$. 

In order to illustrate the difference between the Eulerian and Lagrangian frameworks we consider the residual equation of \Cref{problem::model}. The state equation determining $\eulerState$ is \eqref{eq::residual} with
\begin{equation*}
\begin{aligned}
	R(\Omega_t,\eulerState,v) &= \int_{\Omega_t} \lambda(\pos,\eulerState(\pos)) (\nabla \eulerState \cdot \nabla v) \ddx \\ &- \int_{\Omega_t} f \,v \ddx. 
\end{aligned}
\end{equation*}
A Lagrangian formulation is defined as
\begin{equation}\label{eq::RL}
	R^L(t,\lagrangeState,v^{L,t}) = R(\Omega_t, \lagrangeState\circ T_t^{-1} ,v\circ T_t^{-1} ) .
\end{equation}
Using standard rules of integral transformation and $\lagrangeState\circ T_t^{-1} \circ T_t = \lagrangeState$, \eqref{eq::RL} can be rewritten to
\begin{widetext}
\begin{equation}
	\begin{aligned}	 
		R^L(t,\lagrangeState,v^{L,t}) 
		&= \int_{\Omega} \lambda(T_t(\posX),\lagrangeState) \left[\nabla (\lagrangeState\circ T_t^{-1}) \cdot \nabla (v^{L,t}\circ T_t^{-1})\right]\circ T_t \operatorname{det}\mathbf F_t \ddX \\ & \qquad- \int_{\Omega} f \,v^{L,t} \operatorname{det}\mathbf F_t \ddX.
	\end{aligned}	
\end{equation}
\end{widetext}
\clearpage
Furthermore, using classical transformation rules for gradients (see \Cref{sec::DifferentialOperators}) the fully Lagrangian setting avoiding the occurrence of $T_t^{-1}$ reads:
\begin{widetext}
\begin{equation}
	\begin{aligned}	 
		R^L(t,\lagrangeState,v^{L,t}) &= \int_{\Omega} \lambda(T_t(\posX),\lagrangeState) \left(\nabla \lagrangeState \cdot (\mathbf F_t^{-1}(\mathbf X) \cdot \mathbf F_t^{-\top}(\mathbf X)) \cdot \nabla v^{L,t}\right) \operatorname{det}\mathbf F_t \ddX \\ & \qquad- \int_{\Omega} f \,v^{L,t} \operatorname{det}\mathbf F_t \ddX.
	\end{aligned}	
\end{equation}
\end{widetext}
This "pullback of the shape perturbation to the unperturbed domain" will be detailed for all implemented operators in \Cref{sec::implementation} . 

\newcommand{\someShape}{\omega}

In order to compute the shape derivative, it is customary \cite{allaire2021shape} to introduce a Lagrangian, in an Eulerian setting, by summing the objective function and the state equations
\begin{equation*}
		\mathcal L(\Omega_t,\varphi,\psi) = J(\Omega_t,\varphi) + R(\Omega_t,\varphi,\psi), 
\end{equation*}
where $(\varphi,\psi)\in V(\Omega_t)\times V(\Omega_t)$ are any functions (in the end, $\varphi$ will be replaced by the state $u(\Omega)$ and $\psi$ by the adjoint state). 
The Lagrangian allows us to rewrite the numerator in \eqref{eq::shapeDerivativeDef} as
\begin{align*}
	\shapeFunction(\Omega_t)-\shapeFunction(\Omega) &= \mathcal L(\Omega_t,\eulerState,\psi^{E,t}) - \mathcal L(\Omega,u^0,\psi^0). 
\end{align*}
However, here the drawback is that $(\eulerState,\psi^{E,t})\in V(\Omega_t)\times V(\Omega_t)$  and $(u^0,\psi^0)\in V(\Omega)\times V(\Omega)$, \ie they are not defined over the same functional space. 
It turns out that this Eulerian setting is not easily amenable to automatic differentiation, contrary to the Lagrangian setting that we now introduce. 
Recalling the Lagrangian state \eqref{eq::lagrangeEuler} we define 
the Lagrangian $\mathcal G:\R\times V(\Omega)\times V(\Omega)$ in a Lagrangian setting by
\begin{align}\label{eq::LLagrange}
	\mathcal G(t,\varphi^{L,t}, \psi^{L,t}) &= \mathcal L(T_t(\Omega),\varphi^{L,t}\circ T_t^{-1}, \psi^{L,t}\circ T_t^{-1}) 
\end{align}
Now, we have
\begin{align}\label{eq::shapeDiff}
	\shapeFunction(\Omega_t)-\shapeFunction(\Omega) 
    &= \mathcal G(t,\lagrangeState,\psi^{L,t}) - \mathcal G(0,u^0,\psi^0),
\end{align}
which has the advantage that $(\lagrangeState,\psi^{L,t})$ and ($u^0,\psi^0$) are both defined over the unperturbed domain $\Omega$. Considering \eqref{eq::shapeDiff} in \eqref{eq::shapeDerivativeDef}, and choosing a test function $\psi^{L,t}$ which is independent of $t$, yields by the chain rule
\begin{equation}\label{eq::shapeDerivative2}
	DJ(\Omega)(\mathbf V) = \partial_t\mathcal G + \partial_\phi\mathcal G[\dot{u}^L] ,   
\end{equation}
where
\begin{equation*}
\begin{aligned}
	 \partial_t\mathcal G &= \left(\frac{\partial}{\partial t} \mathcal G(t,u^0,\psi^0)\right)\bigg|_{t=0}, \\ 
  \partial_\phi\mathcal G[\dot{u}^L] &= \left(\frac{d}{d\tau }\mathcal G(0,u^0+\tau\dot{u}^L,\psi^0)\right)\bigg|_{\tau=0}.   
\end{aligned}
\end{equation*}
Here, $\dot{u}^L$ is the Lagrangian shape derivative (also called material derivative) of the state.
Next we introduce the adjoint state $p^0$ with the goal to eliminate the Lagrangian shape derivative of the state. To this end, let $p^0$ be the solution of 
\begin{equation}
\left(\frac{d}{d\tau }\mathcal G(0,u^0+\tau v,p^0)\right)\bigg|_{\tau=0} = 0 \quad \forall v \in V(\Omega).
\end{equation}
Then, for $\psi^0 = p^0$  we have in particular $\partial_\phi\mathcal G[\dot{u}^L]=0$, and \eqref{eq::shapeDerivative2} is reduced to
\begin{equation}\label{eq::derivativeLagrangianSingle}
	DJ(\Omega)(\mathbf V) = \left(\frac{\partial}{\partial t} \mathcal G(t,u^0,p^0)\right)\bigg|_{t=0}. 
\end{equation}
As explained, e.g., in \cite{gangl2021fully} this shape derivative formula is amenable to automatic differentiation. 
\subsection{Shape derivative of multi-physics problems}\label{sec::Multi}
For our multi-physics applications, the objective functionals have the general structure
\begin{equation*}
	\Omega \mapsto \mathcal J(\Omega,u_1(\Omega),\mydots, u_N(\Omega)),
\end{equation*}
where the $N$ scalar or vector-valued fields $u_i$, $i=1,\mydots,N$ are the solutions of the respective governing equations $R_i(\Omega,u_1(\Omega),\mydots, u_N(\Omega),v_i)=0$ for all $v_i \in V_i(\Omega)$. 
For a perturbed domain $\Omega_t$ the perturbed Eulerian states $\eulerState_i \in V_i(\Omega_t)$ satisfy 
\begin{equation}
	R_i(\Omega_t;\eulerState_1,\mydots, \eulerState_N;v_i)=0 \quad \forall v_i \in V_i(\Omega_t).
\end{equation}
In accordance with \eqref{eq::lagrangeEuler} the Lagrangian states $\lagrangeState_i \in V_i(\Omega)$ are defined by 
\begin{equation}\label{eq::lagrangeEulerMulti}
	\lagrangeState_i(\posX) = \eulerState_i \circ T_t(\posX) = \eulerState_i(T_t(\posX)). 
\end{equation}
The Lagrangian is then defined by summing up the objective function and the state equations
\begin{equation*}
	\begin{aligned}
		\mathcal L(\Omega;\varphi_1,\mydots, \varphi_N;&\psi_1,\mydots, \psi_N) = \mathcal J(\Omega;\varphi_1,\mydots, \varphi_N) \\ &+ \sum_{i=1}^N R_i(\Omega;\varphi_1,\mydots, \varphi_N;\psi_i), 
	\end{aligned}
\end{equation*}
for any functions $\varphi_1,\mydots, \varphi_N$ and $\psi_1,\mydots, \psi_N$. 
Recalling \eqref{eq::Hadamard}, $\Omega_t = T_t(\Omega)$, the perturbed Lagrangian in a Lagrangian setting is defined by
\begin{widetext}
\begin{align}\label{eq::LLagrangeMulti}
	\mathcal G(t;\lagrangeStateAny_1,\mydots, \lagrangeStateAny_N; \psi^{L,t}_1,\mydots, \psi^{L,t}_N) = \mathcal L(T_t(\Omega);\lagrangeStateAny_1\circ T_t^{-1},\mydots, \lagrangeStateAny_N\circ T_t^{-1}; \psi^{L,t}_1\circ T_t^{-1},\mydots, \psi^{L,t}_N\circ T_t^{-1}). 
\end{align}    
\end{widetext}
With \eqref{eq::LLagrangeMulti} the analogous arguments from \Cref{sec::Single} can be used to derive the shape derivative formula
\begin{equation}\label{eq::derivativeLagrangianMulti}
	DJ(\Omega)(\mathbf V) = \left(\frac{\partial}{\partial t} \mathcal G(t;u_1^0,\mydots, u_N^0;p_1^0,\mydots, p_N^0)\right)\bigg|_{t=0}, 
\end{equation}
where the adjoint solutions $p_1^0,\mydots, p_N^0$ are determined by
\begin{equation}\label{eq::adjoint}
	\left(\frac{d}{d\tau}\mathcal G(0;u_1^0,\mydots,u_i^0 + \tau \,v_i ,\mydots, u_N^0;p_1^0,\mydots, p_N^0)\right)\bigg|_{\tau=0} = 0, 
\end{equation}
which have to hold for all test functions $v_i\in V_i(\Omega)$. For \eqref{eq::derivativeLagrangianMulti} to hold it is crucial that the adjoint solutions $p^0_1,\mydots, p^0_N$ are determined by \eqref{eq::adjoint} in order to kill terms where the material derivative of the state variables show up.



\section{Implementation}\label{sec::implementation}
In this section, we describe some implementation details of \autofreefem. In view of the theory described in \Cref{sec::Theory}, the symbolic differentiation of expressions plays an important role. 
In particular, differentiation with respect to the perturbation parameter $t$ of the Lagrangian $\mathcal G$ in \eqref{eq::LLagrange}, as well as the Gateaux derivative for the linearization \eqref{eq::Gateaux} and the adjoint problem \eqref{eq::adjoint} have to be performed. Therefore, \autofreefem builds upon the \python package \sympy \cite{sympy}. Beside the symbolic differentiation, the change of variables in \eqref{eq::LLagrange} and the \latex processing uses and extends standard features of \sympy. We give a brief introduction into these topics in \sympy in \Cref{sec::sympy}. In \Cref{sec::fields,sec::TensorAlgebra,sec::DifferentialOperators,sec::Integrals,sec::matrixFunctions,sec::fixedQuantities} we describe the implemented classes.

\subsection{A brief introduction to differentiation, change of variables and \latex processing in \sympy}
\label{sec::sympy}
\begin{figure}[ht]
\begin{subfigure}{0.49\textwidth}
	\tikzstyle{decision} = [diamond, draw, text width=4em, text badly centered, inner sep=0pt]
	\tikzstyle{block} = [rectangle, draw, text width=3em, text centered, rounded corners, text height=0.6em]
	\tikzstyle{line} = [draw, -latex']
	
	\begin{tikzpicture}[node distance = 1.4cm, auto]
		\node [block] (init) {Add};
		\node [block, below left of=init] (decide) {Pow};
		\node [block] (C1) at (-2.3,-2.1) {x};
		\node [block] (C2) at (-0.7,-2.1) {2};
		\node [block, below right of=init] (B2) {Mul};
		\node [block] (C3) at (0.7,-2.1) {x};
		\node [block] (C4) at (2.3,-2.1) {y};
		\path [line] (init) -- (decide);
		\path [line] (decide) -- (C1);
		\path [line] (decide) -- (C2);
		\path [line] (init) -- (B2);
		\path [line] (B2) -- (C3);
		\path [line] (B2) -- (C4);
	\end{tikzpicture}
	\caption{Expression tree for $x^2 + x y$}
	\label{fig::expressionTree}
\end{subfigure}
\begin{subfigure}{0.49\textwidth}
	\tikzstyle{decision} = [diamond, draw, text width=4em, text badly centered, inner sep=0pt]
	\tikzstyle{block} = [rectangle, draw, text width=3.3em, text centered, rounded corners, text height=0.6em]
	\tikzstyle{line} = [draw, -latex']
	
	\begin{tikzpicture}[node distance = 1.5cm, auto]
		\node [block] (init) {Add};
		\node [block, below left of=init] (decide) {Mul};
		\node [block] (C1) at (-2.,-2.1) {2};
		\node [block] (C2) at (-0.5,-2.1) {x};
		\node [block, below right of=init] (B2) {y};
		\path [line] (init) -- (decide);
		\path [line] (decide) -- (C1);
		\path [line] (decide) -- (C2);
		\path [line] (init) -- (B2);
	\end{tikzpicture}
	\caption{Expression tree for $2x+y$}
	\label{fig::expressionTreeD}
\end{subfigure}
\caption{Examples of expression trees}
\label{fig::exptree}
\end{figure}
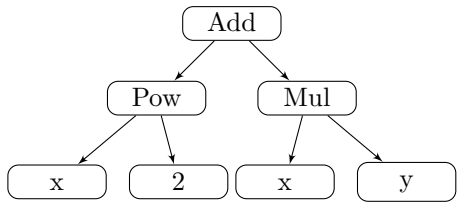
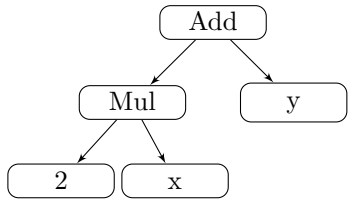

In \sympy, symbolic expressions are stored in expression trees. An expression tree is a data structure with a hierarchical form and the properties:
\begin{enumerate}
	\item Each internal node represents an operator, \eg addition, subtraction, multiplication, division, etc. 
	\item The operands (numbers and variables) are stored in the leaf nodes.
	\item The edges between nodes indicate on which expressions the operators operate.
\end{enumerate}
See \Cref{fig::expressionTree} for a visualization of the expression tree of \textit{expr = x**2 + x*y} in \sympy. 
\paragraph{Differentiation} To differentiate this expression \textit{expr} with respect to the variable $x$ the \sympy command \textit{diff(expr, x)} is used:
\begin{lstlisting}[basicstyle = \ttfamily,columns=fullflexible]
	from sympy import *
	x,y = symbols("x y")
	expr = x**2 + x*y
	dexpr = diff(expr,x)
	print(srepr(expr))
	print(srepr(dexpr))
\end{lstlisting}
Here, the \sympy command \textit{srepr} is used to asses the internal tree representation. The above code gives the outputs:
\begin{quote}
\textit{Add(Pow(Symbol('x'), Integer(2)), Mul(Symbol('x'), Symbol('y')))} \\\\
\textit{Add(Mul(Integer(2), Symbol('x')), Symbol('y'))}
\end{quote}
They correspond to the expression trees in \Cref{fig::exptree} respectively. In order to realize the Gateaux derivative, we use twice the \sympy command \textit{subs} and one time the \sympy command \textit{diff}. The differentiation of \textit{expr} with respect to $x$ into the direction $v$ is given by:
\begin{lstlisting}[basicstyle = \ttfamily,columns=fullflexible]
	v,tau = symbols("v tau")
	expr_vt = expr.subs(x,x+tau*v)
	dexpr_vt = diff(expr_vt,tau)
	dexpr_v = dexpr_vt.subs(tau,0)
	print("expr_vt: ", expr_vt)
	print("dexpr_vt: ", dexpr_vt)
	print("dexpr_v: ", dexpr_v)
\end{lstlisting}
The above code gives the outputs:
\begin{quote}
 \textit{expr\_vt: y*(tau*v + x) + (tau*v + x)**2}\\
 \textit{dexpr\_vt: v*y + 2*v*(tau*v + x)}\\
 \textit{dexpr\_v: 2*v*x + v*y}
\end{quote}
The last output is the sought directional derivative $2 v x + v y$. We remark that the result can be simplified by using the \sympy command \textit{simplify}:
\begin{lstlisting}[basicstyle = \ttfamily,columns=fullflexible]
print("dexpr_v: ", simplify(dexpr_v))
\end{lstlisting}
\begin{quote}
 \textit{dexpr\_v:  v*(2*x + y)}
\end{quote}

\paragraph{Change of variables} In addition to symbolic differentiation, we need to perform a change of variables to obtain the perturbed Lagrangian \eqref{eq::LLagrange} in a Lagrangian framework. To obtain $\mathcal G(t;\lagrangeState_1,\mydots, \lagrangeState_N; p^{L,t}_1,\mydots, p^{L,t}_N)$ in an automatic way, we traverse the expression tree and apply to each operator the corresponding change of variable rule such that $T_t$ and $T_t^{-1}$ cancel out. These rules are are non-trivial transformations for differential operators (see \Cref{sec::DifferentialOperators}) and for integrals (see \Cref{sec::Integrals}).

\paragraph{\latex and \freefem output} The special feature of \autofreefem is that it offers a \latex representation of the input and the derived formulas in coordinate independent direct notation and also a representation of them for the use in \freefem. Remark that \sympy has several built-in options for the output of expressions like the basic string output, a \latex\; output, C code output, and Fortran code output:
\begin{lstlisting}[basicstyle = \ttfamily,columns=fullflexible]
print("String:", expr)
print("Latex:", latex(expr))
print("C code:", ccode(expr))
print("Fortran code:", fcode(expr))
\end{lstlisting}
The above code gives the outputs:
\begin{quote}
\textit{String:  x**2 + x*y} \\
\textit{Latex:  x^{2} + x y} \\
\textit{C code:  pow(x, 2) + x*y} \\
\textit{Fortran code:        x**2 + x*y}
\end{quote}
We are augmenting these built-in output capabilities in two directions. On the one hand, we are extending the \latex processing by introducing coordinate independent direct notation. On the other hand, we introduce the ability to generate code for use in \freefem.

Due to the considerations above, each implemented field and operator is a subclass of the \sympy Function class and thus uses the same mechanisms as elementary functions in \sympy. Additionally, we specify for each class  
\begin{enumerate} 
	\item a rule for the generation of \latex output,  
	\item a rule for the generation of \freefem output, 
	\item if necessary the change of variables to obtain the perturbed Lagrangian functional using the Lagrangian states,
	\item a rule for computing the derivative, 
	\item and if possible some algebraic rules to simplify the expressions.
\end{enumerate}
In the following, we describe each implemented class of \autofreefem. In particular, we give a detailed explanation on how to implement \Cref{problem::model} (see \href{https://gitlab.tugraz.at/autofreefem/autofreefem/-/blob/main/examples/simulation/nonlinearDiffusion.py}{"nonlinearDiffusion.py"} for the full file). To this end, we first import \autofreefem:
\begin{lstlisting}[language=Python,linerange={1-1}]
from autofreefem import *
\end{lstlisting}

\subsection{Fields, Domain, Constants, Expressions and Lagrangian}\label{sec::fields}
In this section, two classes of unknown physical fields, such as temperatures, displacements or velocities, are introduced. These classes are summarized in Table \ref{tab::fields}. 
\begin{table*}[ht]
	\begin{tabular}{ccc}
		\toprule
		description & \autofreefem input & \\
		\midrule
		unknown scalar & ScalarField(symbol, fe-space, mesh, b.c. function, b.c. label) & \\
		unknown vector & VectorField(symbol, fe-space, mesh, b.c. function, b.c. label) & \\
		domain & Domain(\latex symbol, \freefem symbol, boundary 1, boundary 2, \mydots) & \\
		constant & Constant(symbol, numerical value, description, unit) & \\
		expression/abbriviation & Expression(symbol, formula) & \\
		class for computation & Lagrangian(trial fields, test fields, objective, PDE ) & \\
		\bottomrule
	\end{tabular}
	\caption{Implemented classes described in \Cref{sec::fields}}
	\label{tab::fields}
\end{table*}
The \textbf{ScalarField} and the \textbf{VectorField} both take five input arguments. The first argument is a string representing the symbol of the unknown field and is used in the \latex output and the \freefem code. The second argument contains the information how this field should be discretized in \freefem. The third argument is the domain on which the field is defined. The fourth and the fifth argument are related to Dirichlet boundary conditions. In particular, the fourth argument specifies a function for the corresponding values of the Dirichlet boundary data. The fifth argument specifies the boundary labels on which Dirichlet boundary conditions are applied. Thus, for \Cref{problem::model} we make the following definitions:
\begin{lstlisting}[language=Python,linerange={2-7},firstnumber=2]
symbol = 'u'
fespace = 'P1'
mesh = Domain('\\Omega', 'Th')
dirichlet_boundary_function = '0.'
dirichlet_boundary_labels = '4'
u = ScalarField(symbol, fespace, mesh, dirichlet_boundary_function, dirichlet_boundary_labels)
\end{lstlisting}
Here we use conforming finite elements of polynomial degree 1 (P1) for the field $u$. Furthermore, we interoperate homogeneous Dirichlet boundary conditions on all boundaries with label 4. The definitions of the computational domain and of the corresponding mesh are done in the class \textbf{Domain}. The first argument (here '\textbackslash\textbackslash Omega') is the \latex expression, whereas within \freefem code the second argument (here 'Th') will be used. We proceed by specifying the test function $v$:
\begin{lstlisting}[language=Python,linerange={8-8},firstnumber=8]
v = ScalarField('v', fespace, mesh, dirichlet_boundary_function, 'none')
\end{lstlisting}
Note that for the test functions the boundary conditions are inherited form the corresponding unknown fields and therefore the fourth and the fifth argument on line 8 have no effect.

The chosen non-linearity is the diffusion coefficient $\lambda(u) = \nicefrac{\lambda_0}{(1 + u^2)}$: 
\begin{lstlisting}[language=Python,linerange={9-10},firstnumber=9]
Lambda0 = Constant('\\lambda_0', 10., 'diffusion coefficient','')
Lambda = Expression('\\lambda', Lambda0 / (1 + u ** 2))
\end{lstlisting}
In line 9 we first define an object of type \textbf{Constant}. It takes four arguments: a symbol, a numerical value, a description text and a string representing the unit. All constants will be automatically gathered and a \latex table will be generated for the documentation of the used numerical values (see \eg \Cref{tab::FlowParameter} or \Cref{tab::CantileverParameter}). In line 10, we used the class \textbf{Expression}, which has mainly the purpose of introducing an abbreviation to achieve a nicely readable \latex output. It takes two arguments: a symbol and the actual expression.
\begin{remark}
	Note that the classes \textbf{Domain}, \textbf{Constant} and \textbf{Expression} are not necessary in order to set up a simulation in \freefem by \autofreefem. Their purpose is to provide the capability to generate nice \latex output.
\end{remark}
Next we define the bulk source term $f$ using the class \textbf{Constant}:
\begin{lstlisting}[language=Python,linerange={11-11},firstnumber=11]
f = Constant('f', 1., 'source term','')
\end{lstlisting}
For the definition of the variational formulation we rely on the classes \textbf{grad} (see \Cref{sec::DifferentialOperators}), \textbf{inner} (see \Cref{sec::TensorAlgebra}), and \textbf{dx} (see \Cref{sec::Integrals})
\begin{lstlisting}[language=Python,linerange={12-13},firstnumber=12]
a = dx(Lambda*inner(grad(u),grad(v)), mesh)
b = dx(f*v, mesh)
\end{lstlisting}
In order to complete the implementation of \Cref{problem::model}, we set up an object of the class \textbf{Lagrangian} and call the method 'setup\_simulation':
\begin{lstlisting}[language=Python,linerange={14-15},firstnumber=13]
lag = Lagrangian([u], [v], 0, a - b )
lag.setup_simulation('non-linear diffusion')
\end{lstlisting}
The class \textbf{Lagrangian} has four input arguments. The first is a list of all trial fields representing the physical states in the problem, the second input argument is a list of all corresponding test fields. The third argument is the objective function and the last argument is the weak formulation of the PDE constraints.

\subsection{Tensor algebra}\label{sec::TensorAlgebra}
An overview of the three implemented operators of tensor algebra is given in Table \ref{tab::tensorAlgebra}. Let $\{\mathbf e_1,\,\mathbf e_2,\,\mathbf e_3\}$ be the standard Cartesian orthonormal basis. In the present paper, a tensor field $T(\pos) $ of order $k$ assigns to every point $\pos$ a tensor of the form $\underbrace{\R^d\otimes\mydots \otimes \R^d}_{\text{k copies}}$. In this way, we identify scalars as tensors of order zero, vectors as tensors of order one and matrices as tensors of order two. In \autofreefem, the tensor product of two tensors of arbitrary orders $k$ and $k'$, giving rise to a tensor of order $k+k'$, is realized by the class \textbf{TensorProduct}.   
\begin{table}[ht]\centering
	\begin{tabular}{ccc}
		\toprule
		operator & \autofreefem input & \latex \\
		\midrule
		tensor product & TensorProduct(\mydots,~\mydots) & $(\mydots)\otimes(\mydots)$ \\
		dot product &  inner(\mydots,~\mydots) & $(\mydots)\cdot(\mydots)$ \\
		double dot product & inner2(\mydots,~\mydots) & $(\mydots):(\mydots)$ \\
		\bottomrule
	\end{tabular}
	\caption{Implemented operations from tensor algebra. All three operators take two tensor fields as inputs.}
	\label{tab::tensorAlgebra}
\end{table}
Next, we define the dot product (class \textbf{inner}) of two tensors as the contraction of these tensors with respect to the last index of the first one, and the first  index of the second one. For example, the dot product of a third order tensor $\mathbf A = A_{ijl} \mathbf e_i\otimes \mathbf e_j \otimes \mathbf e_l$ and a second order tensor $\mathbf B = B_{km} \mathbf e_k\otimes \mathbf e_m$ gives a third order tensor and reads
\begin{equation*}
	\begin{aligned}
	\mathbf A \cdot \mathbf B &= \left(A_{ijl} \mathbf e_i\otimes \mathbf e_j \otimes \mathbf e_l  \right) \cdot \left( B_{km} \mathbf e_k\otimes \mathbf e_m  \right) \\ 
	&=  A_{ijl}  B_{lm} \mathbf e_i\otimes \mathbf e_j \otimes \mathbf e_m.	
\end{aligned}
\end{equation*}
Here, and in the following, the Einstein summation convention applies. Whenever an index occurs twice, we sum over this index, where Latin indices $i,j,\mydots$ take the values $1,2,3$. 
Furthermore, we define the double dot product (class \textbf{inner2}) of two tensors as the contraction of these tensors with respect to the last two indices of the first one, and the first two indices of the second one. The contraction is performed on the closest indices first, \eg
\begin{equation*}
\begin{aligned}
	\mathbf A: \mathbf B &= \left( A_{ijl} \mathbf e_i\otimes \mathbf e_j \otimes \mathbf e_l \right) : \left(  B_{km} \mathbf e_k\otimes \mathbf e_m \right) \\&=  A_{ijl}  B_{lj} \mathbf e_i.
 \end{aligned}
\end{equation*}
As a consequence of these definitions, we have for second order tensors $\mathbf A,\, \mathbf B,\,\mathbf C$, the relation 
$(\mathbf A \cdot\mathbf B) : \mathbf C = \mathbf A: (\mathbf B \cdot\mathbf C)$. 
The operators in Table \ref{tab::tensorAlgebra} commute with the pull back to the unperturbed domain. 
For the implementation it is also important to note that these operators obey the product rule of differentiation, \ie
\begin{subequations}
\begin{align*}
	(\mathbf A\cdot \mathbf B)' &= \mathbf A' \cdot \mathbf B + \mathbf A \cdot \mathbf B',\\
	(\mathbf A: \mathbf B)' &= \mathbf A': \mathbf B + \mathbf A: \mathbf B', \\
	(\mathbf A \otimes \mathbf B)' &= \mathbf A' \otimes \mathbf B + \mathbf A \otimes \mathbf B',
\end{align*}
\end{subequations}
where $'$ denotes the derivation with respect to a scalar parameter $\tau$.

\subsection{Differential operators}\label{sec::DifferentialOperators}
An overview of the two implemented differential operators is given in Table \ref{tab::DifferentialOperators}.
\begin{table}[ht]
	
	\begin{tabular}{ccc}
		\toprule
		operator & \autofreefem input & \latex \\
		\midrule
		gradient & grad(\mydots) & $\nabla(\mydots)$ \\
		divergence & div(\mydots) & $\operatorname{div}(\mydots)$ \\
		\bottomrule
	\end{tabular}

	\caption{Implemented differential operators}
	\label{tab::DifferentialOperators}
\end{table}
The gradient (class \textbf{grad}) of some scalar-valued function $f:\R^3\rightarrow\R$ is defined as  
\begin{equation}\label{eq::gradient}
	\nabla f(\pos) = 
	\frac{\partial f(\mathbf x)}{\partial x_i} \mathbf e_i
\end{equation}
with the Cartesian coordinates $\pos = (x_1,x_2,x_3)$.
We also use the generalization of the gradient for scalar-valued functions \eqref{eq::gradient} to tensor fields. The gradient of a tensor field $\mathbf A$ of arbitrary order $o$ is defined by
\begin{equation*}
	\nabla\mathbf A(\mathbf x) = \frac{\partial\mathbf A(\mathbf x)}{\partial x_i} \otimes \mathbf e_i.
\end{equation*}
Note that $\nabla\mathbf A$ is a tensor of order $o+1$. For the gradient, the pullback of the shape perturbation can be obtained by application of the chain rule, 
\begin{equation}\label{eq::gradpullback}
	(\nabla_x\mathbf A(\mathbf x))\circ \shapePerturbation(\mathbf X) = \nabla_X\mathbf A(\shapePerturbation(\mathbf X)) \cdot \mathbf F_t^{-1}(\mathbf X).
\end{equation}
The second operator described in this section is the divergence (class \textbf{div}). For a tensor field $\mathbf A$ of order $o\ge 1$, it is given by
\begin{equation}
	\operatorname{div}\mathbf A = 
	\frac{\partial\mathbf A}{\partial x_i} \cdot \mathbf e_i = \nabla\mathbf A(\mathbf x) : \mathbf I 
	.
\end{equation}
Note, that $\operatorname{div}\mathbf A$ is a tensor of order $o-1$ and that the divergence is not defined for a scalar field. 
For the divergence the pull back of the shape perturbation is given by
\begin{equation*}
	\begin{aligned}
	(\operatorname{div}\mathbf A(\pos))\circ \shapePerturbation(\mathbf X) &= (\nabla_x\mathbf A(\mathbf x))\circ \shapePerturbation(\mathbf X) : \mathbf I\\ &= \big(\nabla_X\mathbf A(\shapePerturbation(\mathbf X)) \cdot \mathbf F_t^{-1}(\mathbf X) \big): \mathbf I.	
	\end{aligned}
\end{equation*}

\subsection{Integrals and the normal vector}\label{sec::Integrals}
An overview of the implemented integral operators is given in Table \ref{tab::IntegralOperators}. In \autofreefem domain integrals are understood as integrals over volumes (for 3d problems) or areas (for 2d problems) and are realized by the class \textbf{dx}. This class takes two arguments: the function to be integrated and the integration domain.  
\begin{table}[htb]\centering
	\begin{tabular}{ccc}
		\toprule
		operator & \autofreefem input & \latex \\
		\midrule
		domain integral & dx(f, domain) & $\int_\Omega f \,dx$ \\
		surface/line integral & dsx(f, domain, label) & $\int_\Gamma f \,ds_{x}$ \\
		normal vector & SurfaceNormalVector() & $\mathbf n$ \\
		\bottomrule
	\end{tabular}
	\caption{Implemented integral operators and the normal vector}
	\label{tab::IntegralOperators}
\end{table}
For domain integrals of some tensor field $\mathbf{A}$, the pullback to the unperturbed domain reads
\begin{equation*}
	\int_{\Omega_t} \mathbf A(\pos) \ddx = \int_\Omega (\mathbf A\circ \shapePerturbation(\posX))\; \operatorname{det}\mathbf{F}_t(\posX) \ddX .
\end{equation*}
Boundary integrals (surfaces integrals for 3d problems, line integrals for 2d problems) are realized by the class \textbf{dsx}.
Let $\Gamma$ be part of the boundary of the domain $\Omega$ (characterized by some label). The pullback of an integral over the perturbed boundary $\Gamma_t = T_t(\Gamma)$ of some function $f$ is given by
\begin{equation}\label{eq::boundaryIntegral}
\begin{aligned}
	  \int_{\Gamma_t} f(x) \ddsx = \int_{\Gamma} f(T_t(\posX)) J_\Gamma(\posX) \dd s_X,    
\end{aligned}
\end{equation} 
where the Jacobian determinant is 
\begin{equation*}
    J_\Gamma(\posX) 
    = \operatorname{det}\mathbf F_t(\posX) \;\|\mathbf F_t^{-\top}(\posX) \cdot \mathbf n(\posX)\|,
\end{equation*}
with the normal vector $\mathbf n(\posX)$ to $\Gamma$ at $\posX$. In \autofreefem, the unit exterior normal vector is implemented by the class \textbf{SurfaceNormalVector}.

\subsection{Matrix functions}\label{sec::matrixFunctions}
An overview of the implemented matrix functions is given in Table \ref{tab::MatrixFunctions}. They have in common that they take one matrix, \ie a second order tensor, as input argument.  
\begin{table}[ht]\centering
	\begin{tabular}{ccc}
		\toprule
		operator & \autofreefem input & \latex \\
		\midrule
		matrix transpose & transpose(\mydots) & $(\mydots)^\top$ \\
		matrix trace & tr(\mydots) & $\operatorname{{tr}}\left(\mydots\right)$ \\
		matrix determinant & determinant(\mydots) & $\operatorname{det}\left(\mydots\right)$ \\
		matrix inverse & inverse(\mydots) & $\left(\mydots\right)^{-1}$ \\
		transpose of inverse  & inverse\_transpose(\mydots) & $\left(\mydots\right)^{-\top}$\\
		\bottomrule
	\end{tabular}
	\caption{Implemented matrix functions}
	\label{tab::MatrixFunctions}
\end{table}
The operators in Table \ref{tab::MatrixFunctions} commute with the pull back to the unperturbed domain. 
Furthermore, for the differentiation we have implemented the following rules:
\begin{subequations}
\begin{align*}
	(\mathbf A^\top)' &= (\mathbf A')^\top, \\
	(\operatorname{{tr}}\mathbf A)' &= 	\operatorname{{tr}}(\mathbf A'), \\
	\operatorname{det}(\mathbf A)' &
	= \det \left(\mathbf A \right)  \operatorname{tr} \left (\mathbf A^{-1} \cdot \, \mathbf A'\right), \\
	(\mathbf A^{-1})' &= - \mathbf A^{-1}\cdot \mathbf A' \cdot \mathbf A^{-1}, \\
	(\mathbf A^{-\top})' &= - \mathbf A^{-\top}\cdot (\mathbf A')^\top \cdot \mathbf A^{-\top}.
\end{align*}	
\end{subequations}
\subsection{Fixed quantities}\label{sec::fixedQuantities}
An overview of the fixed quantities implemented in \autofreefem is given in Table \ref{tab::misc}.
\begin{table}[ht]\centering
	\begin{tabular}{ccc}
		\toprule
		quantity & \autofreefem input & \latex \\
		\midrule
		identity matrix & identity() & $\mathbf{I}$  \\
		Cartesian unit vector x-axis  & ex() & $\mathbf e_x$ \\
		Cartesian unit vector y-axis  & ey() & $\mathbf e_y$ \\
		Cartesian unit vector z-axis & ez() & $\mathbf e_z$ \\
		\bottomrule
	\end{tabular}
	\caption{Implemented fixed quantities}
	\label{tab::misc}
\end{table}
The pull back to the unperturbed domain does not alter these quantities and they vanish upon differentiation.
\section{Non-linear and multi-physics shape optimization examples}
\label{sec::examples}
This section presents some examples of shape optimization, solved using the automatic code generation capabilities developed in the present paper. In all examples we use the level-set based mesh evolution method introduced in \cite{allaire2014shape}. For a recent tutorial on this method we refer to \cite{Dapogny2023}.

\subsection{Verification}\label{sec::verification}
For each example, we verified the expressions that we obtained in an automatic way for the shape derivative by looking at the finite difference approximation, as well as the Taylor expansion of the perturbed shape functional. However, the results are shown only for the example of \Cref{sec::Diffusion} in \Cref{fig:convergenceShapeDerivative}. For a fixed shape represented by some chosen level-set function and some chosen fixed vector field $\mathbf V$ we plot the quantities 
\begin{equation*}
	e_1(t) = \left|\frac{J(T_t(\Omega)) - J(\Omega)}{t} - DJ(\Omega)(\mathbf V) \right|,
\end{equation*}
and
\begin{equation}
	e_2(t) = \left|J(T_t(\Omega)) - J(\Omega) - t DJ(\Omega)(\mathbf V) \right|,
\end{equation}
for a sequence of decreasing perturbation parameters $t$. By definition of the shape derivative we have
\begin{equation*}
	e_1(t) = \mathcal{O}(t) \quad \mbox{and} \quad  e_2(t) = \mathcal{O}(t^2) \quad \mbox{as}\; t\searrow0. 
\end{equation*}
We remark that in numerical experiments round-off errors are unavoidable. Thus, for $e_1(t)$ we notice a linear decrease in its magnitude with decreasing $t$ when $t>t^*$, where $t^*$ represents a certain threshold. Conversely, $e_1(t)$ tends to increase for $t<t^*$ due to cancellation errors. For $e_2(t)$ we observe a quadratic decrease rate for decreasing $t$ as long as $t>t^*$ and a more or less constant error measure $e_2$ for $t<t^*$.   

\subsection{Non-linear diffusion}\label{sec::Diffusion}
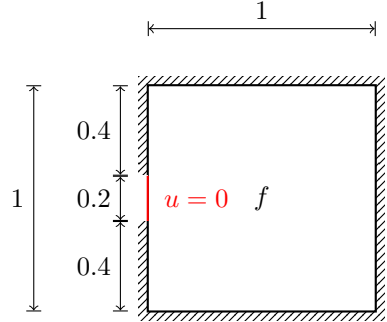
\begin{figure}[t]
\begin{center}
\begin{tikzpicture}[scale=3.]
	\useasboundingbox (-0.7,-0.3) rectangle (1.4,1.5);
	\fill[pattern=north east lines] (-0.04,-0.04) rectangle (1.04,1.04);
	\draw[thick,fill=white] (0,0) rectangle (1,1);
	\fill[white] (-0.06,0.4) rectangle (0,0.6);
	\draw[red] (.03,0.5)  node[right]{$u=0$};
	\draw[red,thick] (0,0.4) --++(0,0.2);
	\draw[|<->|] (0,1.25) --++(1,0) node[pos=0.5,above]{1};
	\draw[|<->|] (-0.5,0.) --++(0.0,1) node[pos=0.5,left]{1};
	\draw[|<->|] (-.12,0.) --++(0.0,0.4) node[pos=0.5,left]{0.4};
	\draw[|<->|] (-.12,0.4) --++(0.0,0.2) node[pos=0.5,left]{0.2};
	\draw[|<->|] (-.12,0.6) --++(0.0,0.4) node[pos=0.5,left]{0.4};
	\node at (0.5,0.5) {$f$};
\end{tikzpicture}	
\end{center}
	\caption{Geometry and boundary conditions of the diffusion problem.}
	\label{fig::diffusion}
\end{figure}
In this first example, we consider an extension of \Cref{problem::model} to a two material shape optimization problem. To this end, we additionally consider a compliance objective function and an area penalization with a fixed Lagrange multiplier $\ell$. The working domain is the unit square, which is heated by a uniform source of magnitude $f$. On a small portion of the left side, we assume Dirichlet boundary conditions, whereas all other boundary parts are assumed to be perfectly isolating (see \Cref{fig::diffusion}). For the non-linear state dependent diffusion coefficient, we assume $\lambda(u) = \chi(x)(1+\alpha u^2)$, where $\alpha$ is a parameter and $\chi(x)$ distinguishes between the two materials. For the material with high conductivity, we have $\chi=1$, whereas for the material with low conductivity, $\chi=0.1$. Note that for $\alpha=0$ the problem becomes linear. A documentation of the formulation of the problem, the adjoint equations and the shape derivative can be obtained by running \href{https://gitlab.tugraz.at/autofreefem/autofreefem/-/blob/main/examples/shapeOptimization/diffusion/run_nonlinearDiffusion.py}{run\_nonlinearDiffusion.py}. The numerical values of the considered physical parameters are supplied in \Cref{tab::nonlineardiffusionParameter}.
\begin{table}[ht]
    \centering
	
\begin{spacing}{1.2}
\begin{tabular}{cccl}
\toprule
Lagrange multiplier & $\ell$ & 100 & $$ \\ 
factor & $\alpha$ & 1.00e-02 & $$ \\ 
source & $f$ & -10 & $$ \\ 
\bottomrule
\end{tabular}
\end{spacing}
		
	\caption{Numerical values of the physical parameters for the non-linear diffusion problem}
	\label{tab::nonlineardiffusionParameter}
\end{table}
\begin{figure}[ht]
	\centering
	\begin{tikzpicture}
		\begin{axis}[width=0.8\linewidth,
            xmode=log,ymode=log,
			grid=both,
			grid style={line width=.1pt, draw=gray!10},
			major grid style={line width=.2pt,draw=gray!50},
			minor tick num=5,
			legend cell align=left,
			legend pos=outer north east,
			legend style={draw=none},
			ylabel=$error$,
			xlabel=$t$]
			\addplot +[] table [x=t,y=numShape,col sep=comma] {convergenceShapenonlineardiffusion.txt};
			\addplot +[] table [x=t,y=taylorShape,col sep=comma] {convergenceShapenonlineardiffusion.txt};
			\addplot  [blue,dashed] coordinates {( 0.1, 0.1*40 ) ( 1e-7, 1e-7*40 ) };
			\addplot  [red,dashed] coordinates {( 0.1, 0.1*0.1 ) ( 1e-7, 1e-7*1e-7 ) }; 
			\legend{$e_1$,$e_2$,$\OO{t}$,$\OO{t^2}$}
		\end{axis}
	\end{tikzpicture}
	\caption{Results of the verification test for the non-linear diffusion example}
	\label{fig:convergenceShapeDerivative}
\end{figure}
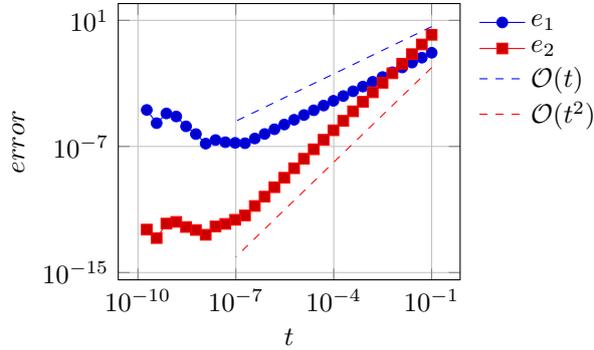
\begin{figure}[ht]
	\begin{subfigure}{0.34\textwidth}
		\centering
		\cropped{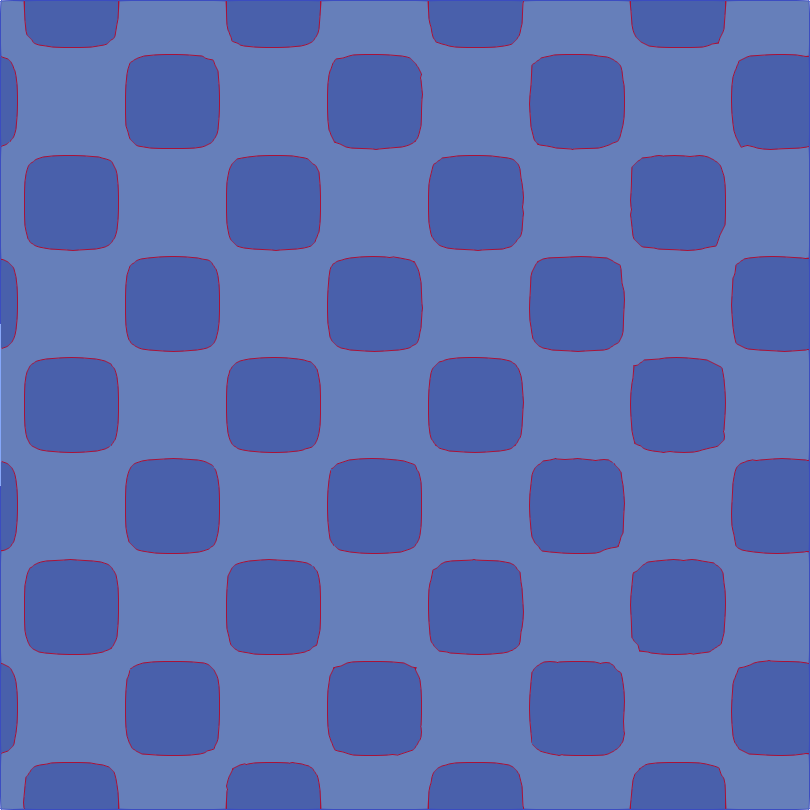}
		\caption{}
		\label{fig:picDiffusionInitialization}
	\end{subfigure}\hfil
	\begin{subfigure}{0.34\textwidth}
		\centering
		\cropped{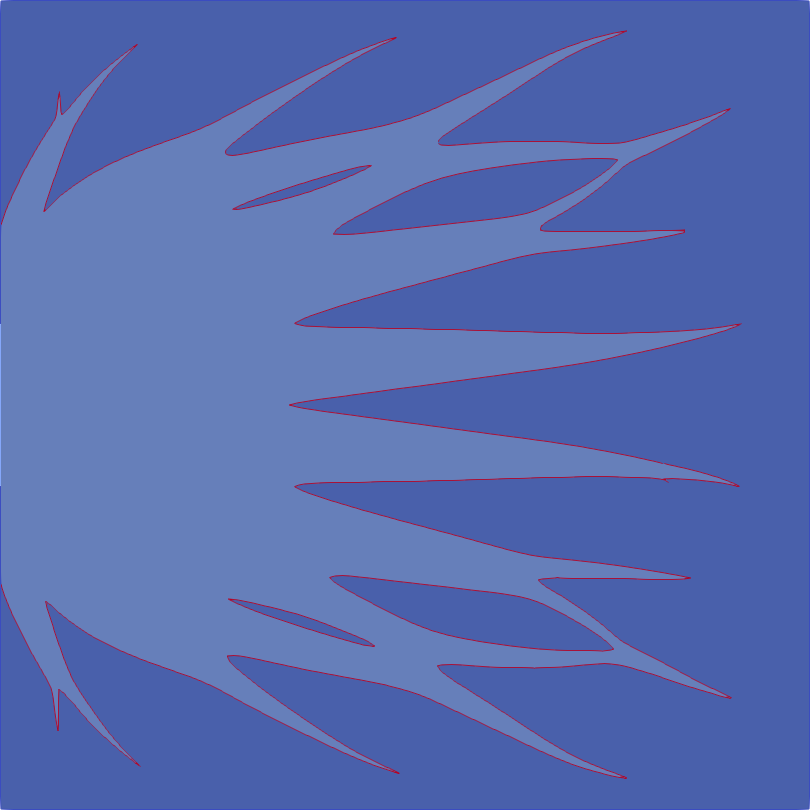}
		\caption{}
		\label{fig:picDiffusionLinearOptimal}
	\end{subfigure}\hfil
	\begin{subfigure}{0.34\textwidth}
		\centering
	\cropped{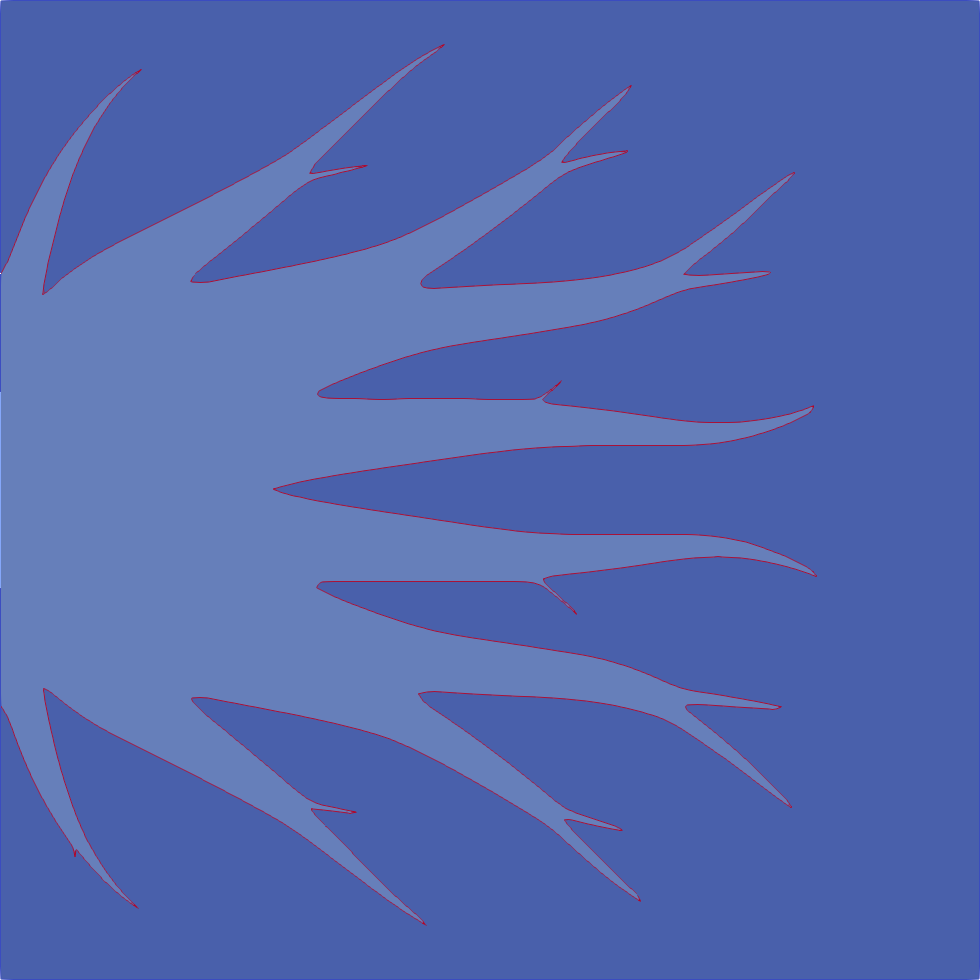}
	\caption{}
	\label{fig:picDiffusionNonLinearOptimal}
	\end{subfigure}
	\caption{Shape optimization of the diffusion problem: (\subref{fig:picDiffusionInitialization}) initialization; (\subref{fig:picDiffusionLinearOptimal}) optimized material distribution for the linear model ($\alpha=0$); (\subref{fig:picDiffusionNonLinearOptimal}) optimized material distribution for the non-linear model ($\alpha=0.01$)}
	\label{fig:picDiffusion}
\end{figure}

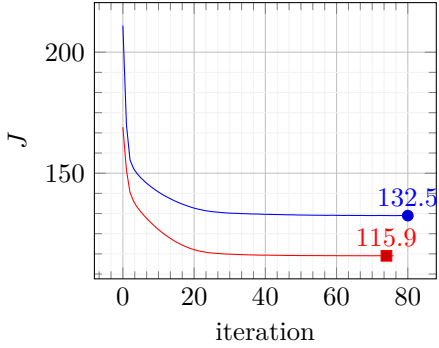
\begin{figure}[ht]
	\centering
	\begin{tikzpicture}
		\begin{axis}[width=0.8\linewidth,
			grid=both,
			grid style={line width=.1pt, draw=gray!10},
			major grid style={line width=.2pt,draw=gray!50},
			minor tick num=5,
			legend cell align=left,
			legend pos=outer north east,
			legend style={draw=none},
			ylabel=$J$,
			xlabel=iteration]
			\addplot +[mark minAbove, every node near coord/.style=] table [x=iteration,y=objective,col sep=comma] {objectivelineardiffusion.txt};
			\addplot +[mark minAbove, every node near coord/.style=] table [x=iteration,y=objective,col sep=comma] {objectivenonlineardiffusion.txt};
		\end{axis}
	\end{tikzpicture}
	\caption{Convergence history for the diffusion problem: The blue curve shows the evolution of the objective function for the linear model, whereas the red curve for the non-linear model.}
	\label{fig:objectivenonlinearScalarLaplace}
\end{figure}

The initialization and the corresponding computed optimal designs for a linear model ($\alpha=0$) and the described non-linear model are displayed in \Cref{fig:picDiffusion}. The evolution of the objective function is reported in \Cref{fig:objectivenonlinearScalarLaplace}. We observe that for the non-linear model the obtained value of the objective function is lower than for the linear model. This was  expected because, in the nonlinear model, the diffusion coefficient is larger than in the linear model. 

\subsection{Non-linear Elasticity}\label{sec::Elasticity}
In this section, we revisit the elasticity cantilever problem discussed in \Cref{sec::linearCantilever}, but now explore both geometrically non-linear and material non-linear behaviors.
Again, the working domain is a rectangle of size $2 \times 1$, with zero displacement boundary condition on the left side and a vertical load applied on a small portion of $0.1$ at the middle of the right side denoted by $\Gamma_N$ such that the resultant force has unit magnitude. All other sides are traction free. The geometry and the boundary conditions are illustrated in \Cref{fig::Cantilever}. There are no body forces.
The objective function is analogously to \Cref{sec::linearCantilever} (see also \Cref{appendix::Cantilever}) the sum of compliance and a fixed Lagrange multiplier $\ell$ multiplied by the area of the solid,
\begin{align*}
J( {\mathbf{ u } }) &=  \int_{D} \ell \raisebox{\depth}{$\chi$} \,dx +  \int_{D} \mathbf S( {\mathbf{ u } }) \mathrel{:} \mathbf E( {\mathbf{ u } })  \,dx,
 \end{align*}
 with the second Piola-Kirchhoff stress tensor $\mathbf S$ and the Green-Lagrange strain tensor $\mathbf E$.

\subsubsection{Non-linear elasticity with Saint Venant-Kirchhoff material}
In this section, we consider geometrically non-linear elasticity with the (linear) Saint Venant-Kirchhoff material (see also \cite[Section 8]{allaire2004structural}). The formulation of the non-linear elasticity problem, the adjoint equations and the shape derivative can be obtained by running \href{https://gitlab.tugraz.at/autofreefem/autofreefem/-/blob/main/examples/shapeOptimization/elasticity/run_nonlinearElasticitySaintVernant.py}{run\_nonlinearElasticitySaintVernant.py}. The numerical values of the considered physical parameters are supplied in \Cref{tab::CantileverParameter}. 

\begin{figure}[ht]
	\centering
	\begin{subfigure}{0.49\textwidth}
		\cropped{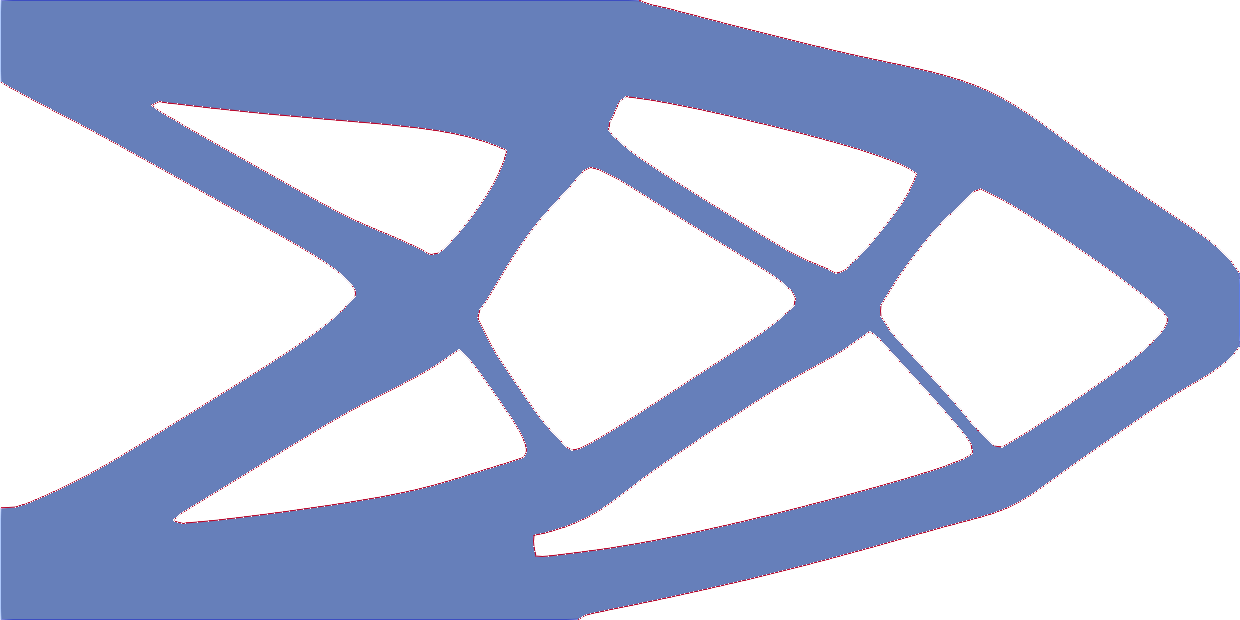}
	\end{subfigure}
	\caption{Optimal design for the non-linear elastic cantilever with Saint Venant-Kirchhoff material.}
	\label{fig:picElasticSaint}
\end{figure}

\begin{figure}[ht]
	\centering
	\begin{tikzpicture}
		\begin{axis}[width=0.9\linewidth,
			grid=both,
			grid style={line width=.1pt, draw=gray!10},
			major grid style={line width=.2pt,draw=gray!50},
			minor tick num=5,
			legend cell align=left,
			legend pos=north east,
			legend style={draw=none},
			ylabel=$J$,
			xlabel=iteration,ymin=0.45,ymax=0.7]
			\addplot +[mark endAbove, every node near coord/.style=] table [x=iteration,y=objective,col sep=comma] {objectivelinearElasticity.txt};
			\addplot+[mark endBelow, every node near coord/.style=] table [x=iteration,y=objective,col sep=comma] {objectivenonlinearElasticitySaintVernant.txt};
			\addplot+[mark end, every node near coord/.style=] table [x=iteration,y=objective,col sep=comma] {objectivenonlinearElasticityNeoHookean.txt};
			\legend{linear,Saint-Venant,Neo-Hookean}
		\end{axis}
	\end{tikzpicture}
	\caption{Convergence history for the cantilever problem with different material laws}
	\label{fig:convergenceElasticity}
\end{figure}
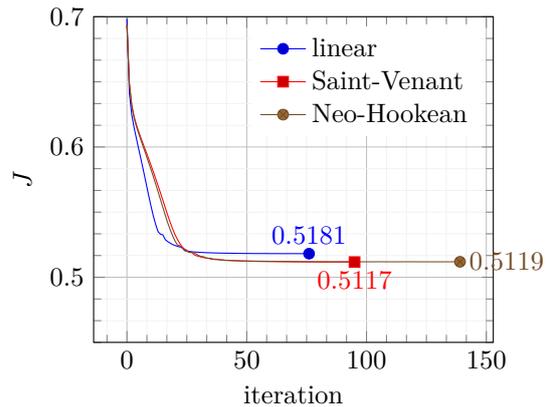

As initialization we use the same geometry as for the linear elastic case (see \Cref{fig:picCantileverElasticInitialization}). The computed optimal design is displayed in \Cref{fig:picElasticSaint}. Due to the non-linear model the design is not symmetric with respect to a horizontal line as it was for the linear model. The evolution of the objective function is reported in \Cref{fig:convergenceElasticity}.

\subsubsection{Non-linear elasticity with Neo-Hookean material}
In this section, we consider now a geometrically and materially non-linear elasticity formulation by resorting to a Neo-Hookean material law \cite{wriggers2008nonlinear}. The formulation of the non-linear elasticity problem, the adjoint equations and the shape derivative can be obtained by running \href{https://gitlab.tugraz.at/autofreefem/autofreefem/-/blob/main/examples/shapeOptimization/elasticity/run_nonlinearElasticityNeoHookean.py}{run\_nonlinearElasticityNeoHookean.py}. The numerical values of the considered physical parameters are supplied in \Cref{tab::CantileverParameter}.

\begin{figure}[ht]
	\begin{subfigure}{0.49\textwidth}
		\cropped{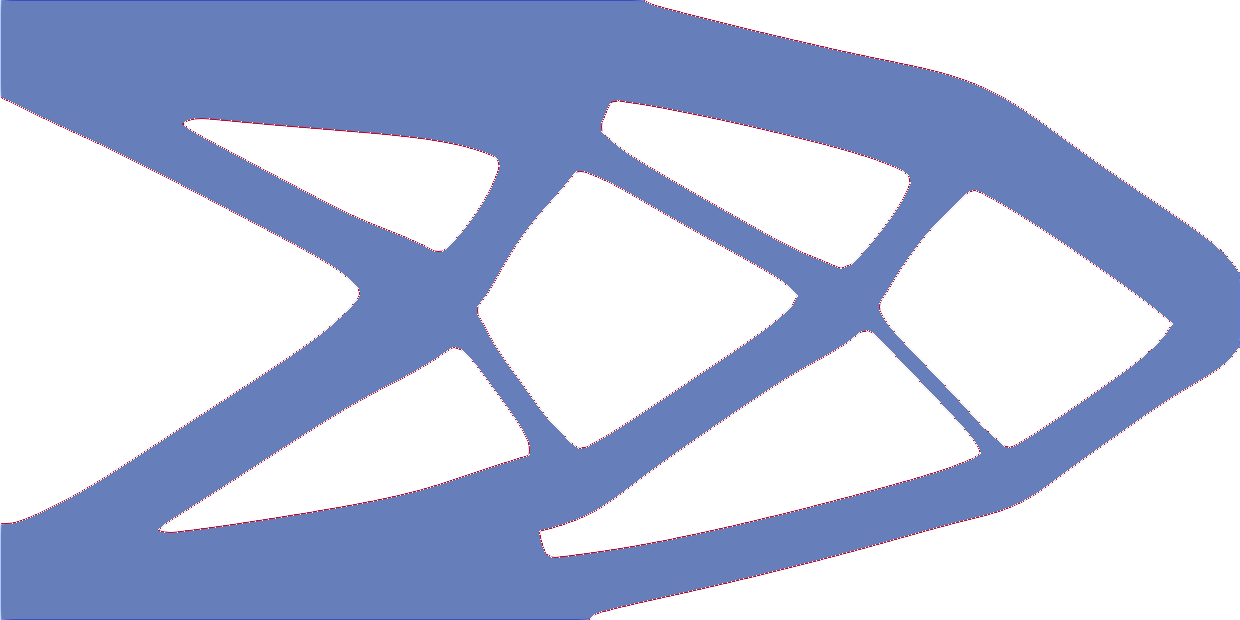}
		\caption{}
		\label{fig:picElasticNeoMin}
	\end{subfigure}
	\begin{subfigure}{0.49\textwidth}
		\cropped{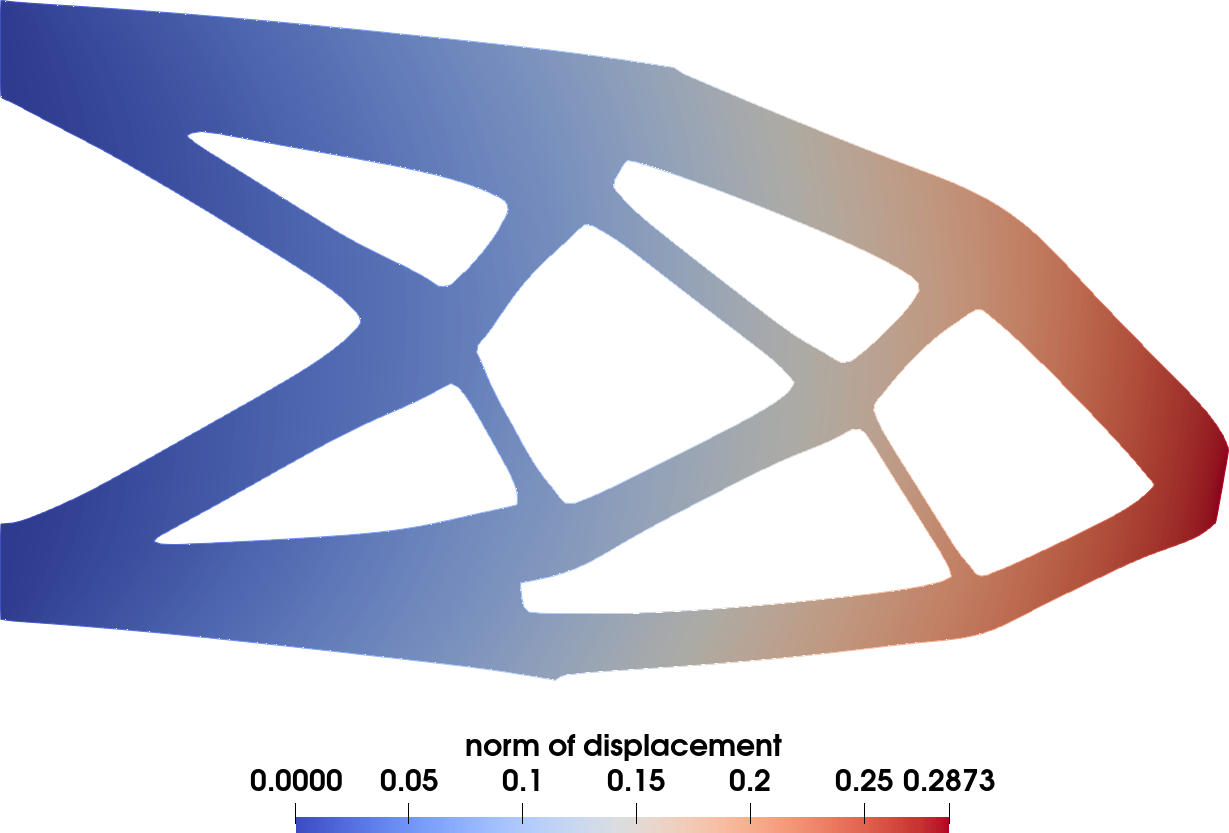}
		\caption{}
		\label{fig:picElasticNeoDeformed}
	\end{subfigure}
	\caption{Non-linear elastic cantilever with Neo-Hookean material: (\subref{fig:picElasticNeoMin}) undeformed optimal design; (\subref{fig:picElasticNeoDeformed}) deformed optimal design. The colors indicate the norm of the dispacement.}
	\label{fig:picElasticNeo}
\end{figure}

As initialization we use the same geometry as for the linear elastic case (see \Cref{fig:picCantileverElasticInitialization}). The computed optimal design is displayed in \Cref{fig:picElasticNeoMin}. Again, due to the non-linear model, the design is not symmetric with respect to a horizontal line as it was for the linear model. Furthermore, the optimal design differs from the optimal design obtained for the Saint Venant-Kirchhoff material law. In \Cref{fig:picElasticNeoDeformed} the deformed optimal design is shown. The evolution of the objective function is reported in \Cref{fig:convergenceElasticity}. Although the optimal designs for the different models differ from each other, the obtained values of the objective functions are quite similar.
\subsection{Thermo-elasticity}\label{sec::ThermoElasticity}
\begin{figure*}[htb]
	\centering
	\begin{tikzpicture}[scale=1.7]
			\useasboundingbox (-0.7,-0.3) rectangle (4.4,1.5);
			\fill[pattern=north east lines] (-0.1,0) rectangle (4.1,1);
			\draw[thick,fill=gray!20] (0,0) rectangle (4,1);
			\node at (0.2,0.5) {$\Gamma_D$};
			\node at (3.8,0.5) {$\Gamma_D$};
			\node at (2,0.4) {$\rho$, $\nu$, $E$};		
			\draw[red,thick] (0,1.0) --++(4,.0) node[pos=0.5,below]{$\Gamma_N$ \quad$T=0^\degree/30^\degree/-30^\degree$};
			\draw[blue,thick] (0,0.0) --++(4,.0) node[pos=0.5,below]{$\Gamma_1$ \quad$T=0^\degree/0^\degree/-30^\degree$};
			\draw[|<->|] (0,1.7) --++(4,0) node[pos=0.5,above]{$4m$};
			\draw[|<->|] (-0.5,0.) --++(0.0,1) node[pos=0.5,left]{$1m$};
			\draw[-latex] (4.3,1) -- ++(0,-0.4)node[right]{$g$};
			\foreach \x in {0,0.4,...,4} {
				\draw[-latex] (\x,1.3) -- ++(0,-0.3);
			}
			\draw (0,1.3) -- (4,1.3) node[midway,above]{$q$};
	\end{tikzpicture}
	\caption{Geometry and boundary conditions of the thermo-elastic bridge problem. }
	\label{fig::thermoElasticity}
\end{figure*}
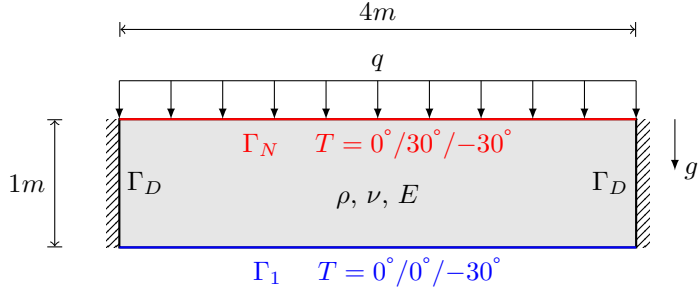
In this example, we optimize a bridge, which is mechanically loaded as well as experiences deformations due to a temperature change. The geometry and the boundary conditions are illustrated in \Cref{fig::thermoElasticity}. The working domain $D$ is a rectangle of size $4m \times 1m$, with zero displacement boundary condition on the left and right sides $\Gamma_D$. A vertical load of constant magnitude $q$ is applied on the top edge $\Gamma_N$ of the domain. Furthermore, the self-weight (density $\rho$, gravitational acceleration $g$) of the bridge is taken into account. The bottom side is traction free. 
\begin{table}[ht]
	\begin{center}
		
\begin{spacing}{1.2}
\begin{tabular}{cccl}
\toprule
Poisson's ratio & $\nu$ & 0.23 & $-$ \\ 
Young's modulus & $E_{0}$ & 3.20e+07 & $\frac{N}{m^2}$ \\ 
density & $\rho$ & 2.50e+03 & $\frac{kg}{m^3}$ \\ 
disp. cost factor & $\gamma$ & 100 & $\frac{N}{m^2}$ \\ 
gravitational acc. & $g$ & 9.81 & $\frac{m}{s^2}$ \\ 
material cost factor & $\ell$ & 1 & $\frac{N}{m^2}$ \\ 
penalty parameter & $\beta_{T}$ & 1.00e+11 & $\frac{W}{K m^2}$ \\ 
penalty parameter & $\beta$ & 1.00e+11 & $\frac{N}{m^3}$ \\ 
ther. conductivity & $k_{0}$ & 1.25 & $\frac{W}{K m}$ \\ 
ther. expansion coeff. & $\alpha_{T}$ & 3.60e-05 & $\frac{1}{K}$ \\ 
vert. load comp. & $q$ & 5.00e+03 & $\frac{N}{m^2}$ \\ 
\bottomrule
\end{tabular}
\end{spacing}

	\end{center}
	\caption{Numerical values of the physical and numerical parameters for the thermo-elasticity problem}
	\label{tab::thermoElasticityParameter}
\end{table}
For the thermal part of the problem we prescribe the temperature change $T$ on the lower and the upper edges $\Gamma_T = \Gamma_N \cup \Gamma_1$ and consider three different cases: (a) no temperature change ($g_{T}(x)=0$), (b) $T=30^\degree$ on the upper edge and $T=0^\degree$ on the lower edge ($g_{T}(x)=30y$), and (c) $T=-30^\degree$ on both edges ($g_{T}(x)=-30$). The state $(T,\, \mathbf u) \in H^1(D)\times [H^1(D)]^2$ is the solution of the classical one-sided coupled thermo-elasticity problem \cite{NowackiWitold1986}
\begin{align}
 &\int_{D} k \left(  \nabla \delta T  \cdot  \nabla T  \right) \,dx +\int_{\Gamma_{T}} \beta_{T} (T-g_{T}(x)) \delta T \,ds_x = 0, \\
 &\int_{\Gamma_D} \beta \left( \delta{\mathbf{ u}} \cdot  {\mathbf{ u } } \right) \,ds_x  +  \int_{D} \pmb \sigma( {\mathbf{ u } },T) \mathrel{:}  \pmb\epsilon(\delta{\mathbf{ u}}) \,dx \nonumber\\ &=   \int_{\Gamma_N} q \left( \mathbf e_y \cdot \delta{\mathbf{ u}} \right) \,ds_x +  \int_{D} \rho g \raisebox{\depth}{$\chi$} \left( \mathbf e_y \cdot \delta{\mathbf{ u}} \right) \,dx,
\end{align}
 for all test functions $(\delta T,\, \delta{\mathbf{ u}}) \in H^1(D)\times [H^1(D)]^2$. We have used the following abbreviations
 \begin{align*}
\pmb\epsilon( {\mathbf{ u } }) &= \frac{ \nabla  {\mathbf{ u } }  +   \nabla  {\mathbf{ u } } ^\top }{2},\\
\pmb\epsilon_{e}( {\mathbf{ u } },T) &=  \pmb\epsilon( {\mathbf{ u } }) - \alpha_{T} T \mathbf{I},\\
\pmb \sigma( {\mathbf{ u } },T) &= \lambda \mathbf{I}  \,\operatorname{tr}\left(\pmb\epsilon_{e}( {\mathbf{ u } },T)\right)  + 2 \mu \pmb\epsilon_{e}( {\mathbf{ u } },T),
 \end{align*}
 where $\alpha_{T}$ is the isotropic thermal expansion coefficient.
 For the domain occupied by material we have $k = k_0$, $E = E_0$, and $\raisebox{\depth}{$\chi$}=1$. Contrary, for the void space we have assumed $k = k_0/1000$, $E = E_0/1000$, and $\raisebox{\depth}{$\chi$}=0$. 
 Note that for the imposition of Dirichlet boundary conditions the penalty method is used (penalty parameter $\beta$ for $\mathbf u = \mathbf 0$ on $\Gamma_D$, and $\beta_T$ for $T=g_T(x)$ on $\Gamma_T$).
 This allows to easily post-process the bearing forces for the evaluation of the objective function.

The numerical values of the considered physical and numerical parameters for concrete material are supplied in \Cref{tab::thermoElasticityParameter}.
The objective is to minimize the following three effects:
\begin{itemize}
	\item the vertical deformation of the upper edge of the domain (displacement cost factor $\gamma$),
	\item the horizontal bearing forces on $\Gamma_D$,
	\item and the material consumption measured as the area (material cost factor $\ell$). 
\end{itemize}
The precise objective function to be minimized is 
\begin{align*}
J( {\mathbf{ u } }) &=  -  \int_{\Gamma_N} \gamma \left( \mathbf e_y \cdot  {\mathbf{ u } } \right) \,ds_x \\&+\int_{\Gamma_D} \beta \left( \mathbf e_x \cdot  {\mathbf{ u } } \right)^2 \,ds_x  +  \int_{\Omega} \ell \raisebox{\depth}{$\chi$} \,dx.
 \end{align*}
The the adjoint equations and the shape derivative can be obtained by running \href{https://gitlab.tugraz.at/autofreefem/autofreefem/-/blob/main/examples/shapeOptimization/thermoElasticity/run_thermoElastic.py}{run\_thermoElastic.py}. 
\begin{figure}[ht]
	\begin{subfigure}{0.91\linewidth}\centering
		\cropped{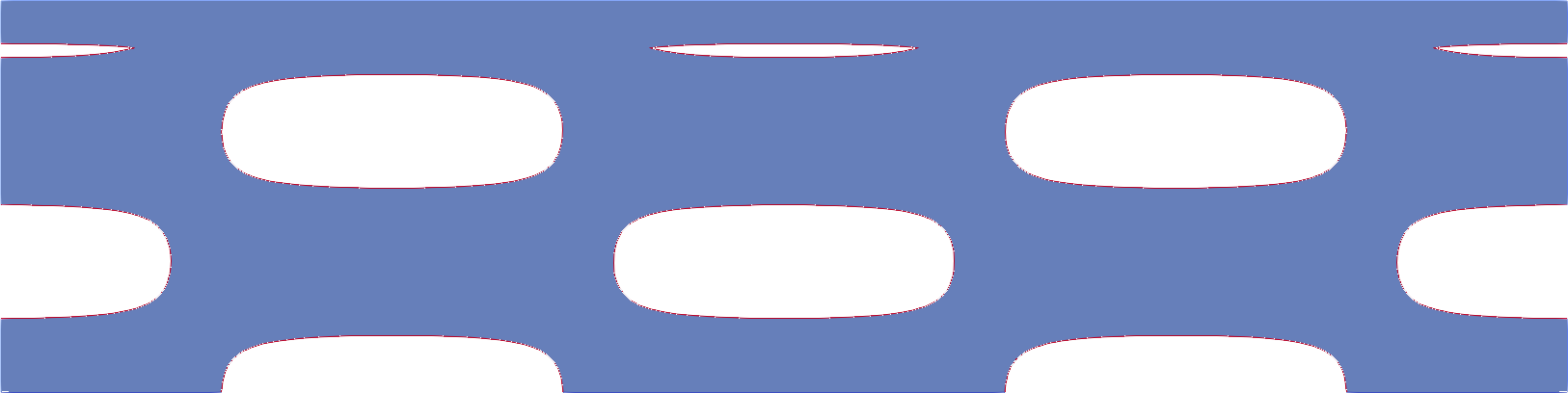}
		\caption{}
		\label{fig:picthermoelasticityInitial}
	\end{subfigure}
	\begin{subfigure}{0.92\linewidth}\centering
		\cropped{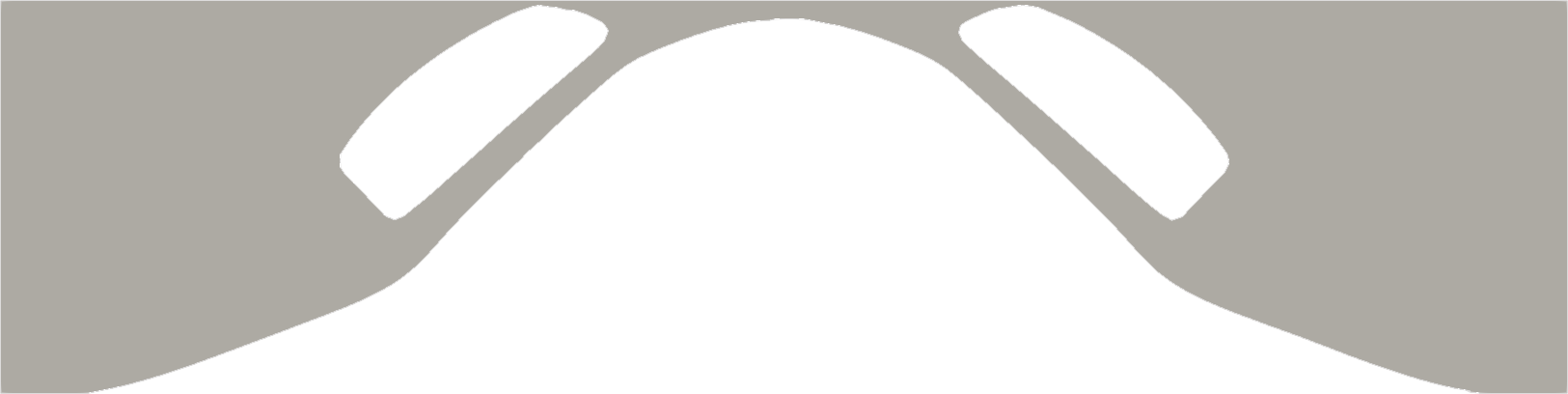}
		\caption{}
		\label{fig:picthermoelasticityA}
	\end{subfigure}
	\begin{subfigure}{0.99\linewidth}\centering
		\cropped{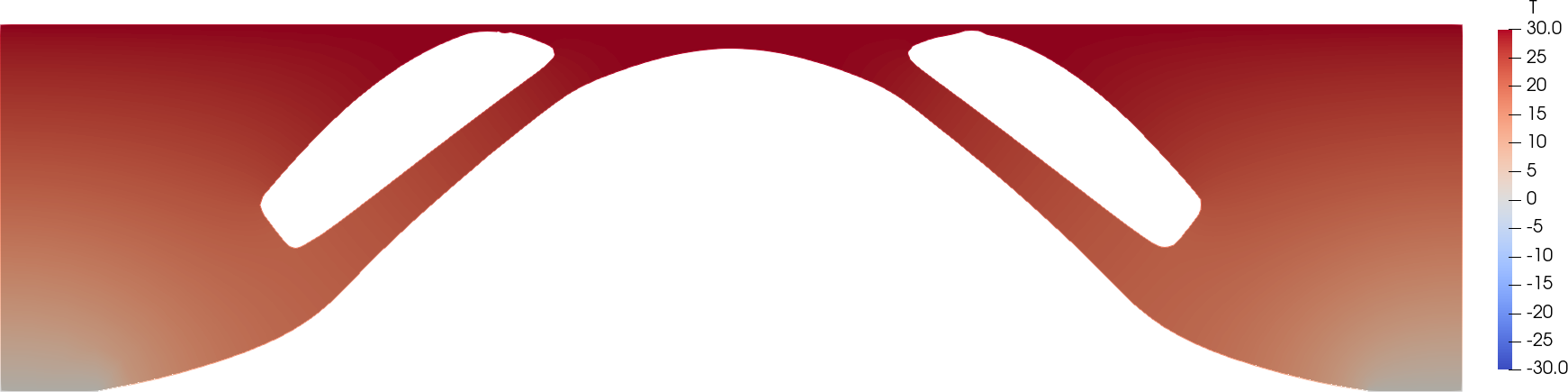}
		\caption{}
		\label{fig:picthermoelasticityB}
	\end{subfigure}
	\begin{subfigure}{0.92\linewidth}\centering
		\cropped{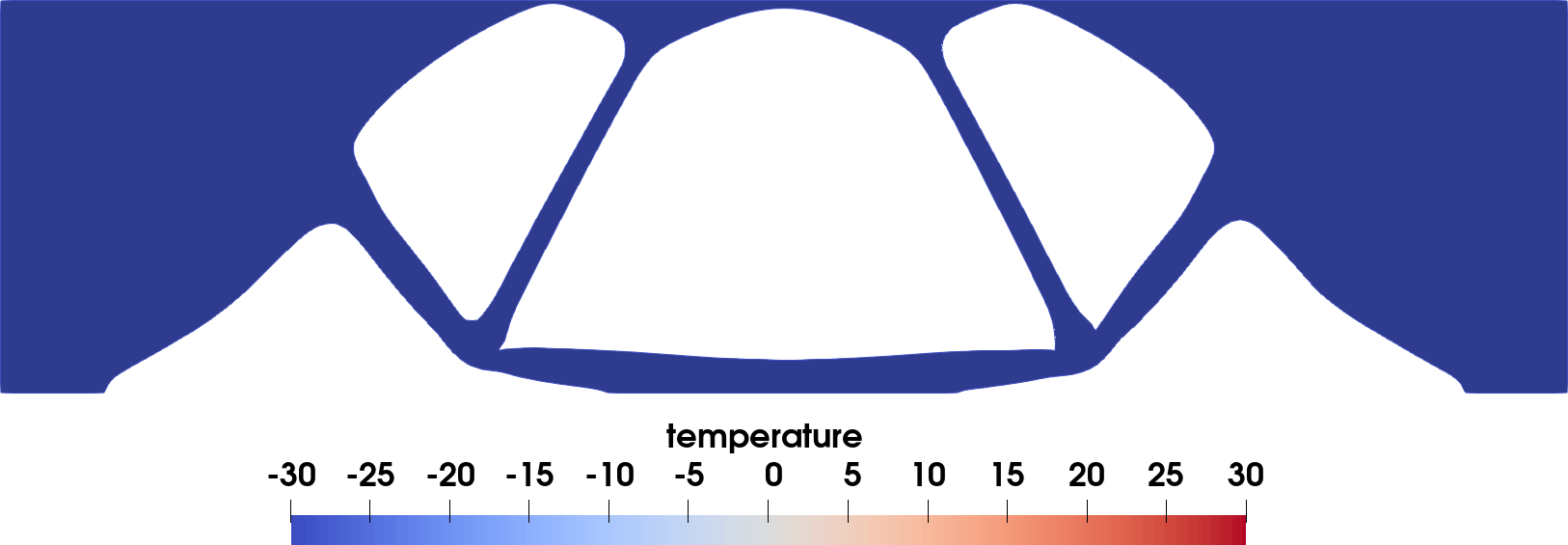}
		\caption{}
		\label{fig:picthermoelasticityC}
	\end{subfigure}
	\caption{(\subref{fig:picthermoelasticityInitial}) initialization; (\subref{fig:picthermoelasticityA}) optimized design for $T=0\degree$; (\subref{fig:picthermoelasticityB}) optimized design for $T=30\degree$; (\subref{fig:picthermoelasticityC}) optimized design for $T=-30\degree$. The colors in optimized designs indicate the temperature distribution over the structure.}
	\label{fig:picthermoelasticity}
\end{figure}

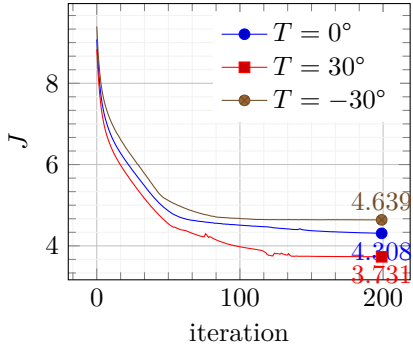
\begin{figure}[ht]
	\centering
	\begin{tikzpicture}
		\begin{axis}[width=0.8\linewidth,
			grid=both,
			grid style={line width=.1pt, draw=gray!10},
			major grid style={line width=.2pt,draw=gray!50},
			minor tick num=5,
			legend cell align=left,
			legend pos=north east,
			legend style={draw=none},
			ylabel=$J$,
			xlabel=iteration]
			\addplot +[mark endBelow, every node near coord/.style=] table [x=iteration,y=objective,col sep=comma] {objectivethermoelastic.txt};
			\addplot +[mark endBelow, every node near coord/.style=] table [x=iteration,y=objective,col sep=comma] {objectivethermoelastic30b.txt};
			\addplot +[mark endAbove, every node near coord/.style=] table [x=iteration,y=objective,col sep=comma] {objectivethermoelastic30c.txt};
			\legend{$T=0\degree$,$T=30\degree$,$T=-30\degree$}
		\end{axis}
	\end{tikzpicture}
	\caption{Convergence history for the thermo-elasticity problem}
	\label{fig:convergenceThermoElasticity}
\end{figure}

The chosen initial design is depicted in \Cref{fig:picthermoelasticityInitial}. The optimized designs for the three load cases are visualized in \Cref{fig:picthermoelasticityA,fig:picthermoelasticityB,fig:picthermoelasticityC}. The shapes obtained for load cases (a) and (b) exhibit remarkable similarity, while load case (c) yields a significantly different design. In the latter scenario, the optimal design retains the lower horizontal part in the middle of the domain. This phenomenon is due to temperature shrinkage, which creates an uplift force in this particular region counteracting the loading $q$ and the dead load. The evolution of the objective function for the three load cases are reported on \Cref{fig:convergenceThermoElasticity}. We note that the rise in temperature in load case (b) positively impacts the objective function, in contrast to the temperature decrease in load case (c). 
\subsection{Fluid-structure interaction}\label{sec::fluidStructureExample}
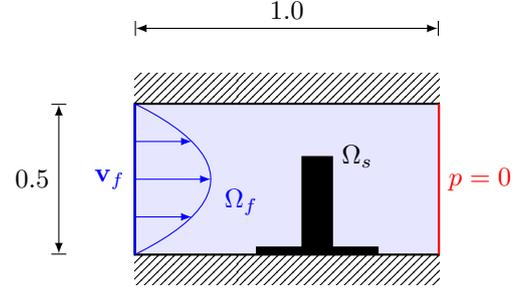
\begin{figure}[htb]
	\centering
	\begin{tikzpicture}[scale=4.,>=latex]
		\def\h{0.5}
		\def\hfound{0.025}
		\def\hstruct{0.3}
		\def\bleft{0.4}
		\def\bb{0.15}
		\def\b{0.1}
		\def\ll{1.0}
		\useasboundingbox (-0.7,-0.3) rectangle (3.4,1.5);
		\fill[pattern=north east lines] (0.,-0.1) rectangle (\ll,\h+0.1);
		\draw[thick,fill=blue!10] (0,0) rectangle (\ll,\h) node[pos=0.35,blue]{$\Omega_f$};		
		\draw[red,thick] (\ll,0.0) --++(0,\h) node[pos=0.5,right]{$p=0$};
		
		\draw[blue,thick] (0,0.0) --++(0,\h) node[pos=0.5,left]{$\mathbf v_f$};
		\draw[blue,domain=0:1] plot (\x-\x*\x,\h*\x) ;
		\foreach\x in {0.25,0.5,0.75}{\draw[blue,->] (0,\h*\x) -- (\x-\x*\x,\h*\x);};
		\draw[|<->|] (0,\h+.25) --++(\ll,0) node[pos=0.5,above]{\ll};
		\draw[|<->|] (-0.25,0.) --++(0.0,\h) node[pos=0.5,left]{\h};
		\draw[fill] (\bleft,0.) -- ++ (0.,\hfound) -- ++(\bb,0) -- ++ (0.,\hstruct) -- ++(\b,0) node[right]{$\Omega_s$} -- ++ (0.,-\hstruct)  -- ++(\bb,0)-- ++ (0.,-\hfound) -- cycle;
	\end{tikzpicture}
	\caption{Geometry and boundary conditions of the fluid-structure interaction problem. }
	\label{fig::fluidstructure}
\end{figure}

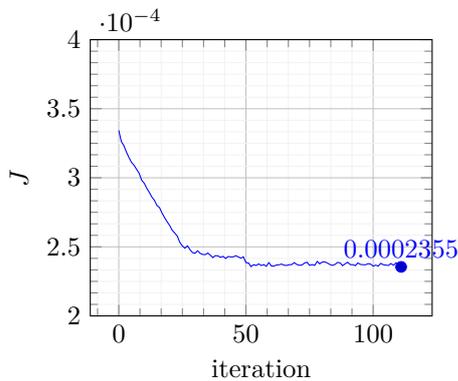
\begin{figure}[ht]
	\centering
	\begin{tikzpicture}
		\begin{axis}[width=0.8\linewidth,
			grid=both,
			grid style={line width=.1pt, draw=gray!10},
			major grid style={line width=.2pt,draw=gray!50},
			minor tick num=5,
			legend cell align=left,
			legend pos=outer north east,
			legend style={draw=none},
			ylabel=$J$,
			xlabel=iteration,ymin=0.2e-3,ymax=0.4e-3]
			\addplot+[mark minAbove, every node near coord/.style=] table [x=iteration,y=objective,col sep=comma] {objectivefluidstructureinteractionnonlinear.txt};
		\end{axis}
	\end{tikzpicture}
	\caption{Convergence history for the fluid-structure problem}
	\label{fig:convergenceFluidStructure}
\end{figure}

In this section, we describe a fluid-structure interaction example which is motivated by \cite{Yoon2010,YOON2014499} and \cite{feppon2019shape}. Here, we assume non-linear fluid flow and non-linear elastic deformations by the structure. We use an arbitrary Lagrangian-Euler formulation (ALE) \cite{le2001fluid,dowell2001modeling,hou2012numerical} and therefore four unknown fields are sought: the elastic displacement field  $\mathbf u$, the fluid velocity $\mathbf v$, the fluid pressure $p$, and an extension of the displacement field to the fluid domain $\mathbf u_{ext}$. The geometry and boundary conditions of the problem are illustrated in \Cref{fig::fluidstructure}.
The formulation of the fluid-structure interaction problem, the adjoint equations and the shape derivative can be obtained by running \href{https://gitlab.tugraz.at/autofreefem/autofreefem/-/blob/main/examples/shapeOptimization/fluidStructureInteraction/run_FluidStructureInteractionNonlinear.py}{run\_FluidStructureInteractionNonlinear.py}. Note that the problem is quite complicated and we give only the objective function
\begin{align}
           J = \int_{\Omega_{s}} \ell \,dx + \int_{\Omega_{s}}  \mathbf S( {\mathbf{ u } }) \mathrel{:} \mathbf E( {\mathbf{ u } }) \,dx,
\end{align}
with the second Piola-Kirchhoff stress tensor $\mathbf S$ and the Green-Lagrange strain tensor $\mathbf E$.
The numerical values of the considered physical parameters are supplied in \Cref{tab::fluidStructureInteractionParameter}.
\begin{table}[ht] 
		\begin{center}
		
\begin{spacing}{1.2}
\begin{tabular}{cccl}
\toprule
Lagrange multiplier & $\ell$ & 5.00e-03 & $\frac{N}{m^2}$ \\ 
Lamé constant & $\lambda$ & 0.2645 & $\frac{N}{m^2}$ \\ 
Lamé constant & $\mu$ & 2.38 & $\frac{N}{m^2}$ \\ 
coupling parameter & $\gamma_{f}$ & 1.00e+08 & $\frac{N}{m^3}$ \\ 
fluid density & $\rho$ & 1 & $\frac{kg}{m^3}$ \\ 
fluid viscosity & $\mu_{T}$ & 5.00e-03 & $Pa s$ \\ 
penalty parameter & $\epsilon$ & 1.00e-08 & $\frac{1}{Pa s}$ \\ 
\bottomrule
\end{tabular}
\end{spacing}

	\end{center}
	\caption{Numerical values of the physical parameters for the fluid-structure interaction problem}
	\label{tab::fluidStructureInteractionParameter}
\end{table}
\begin{figure*}
	\begin{subfigure}{0.49\textwidth}
		\cropped{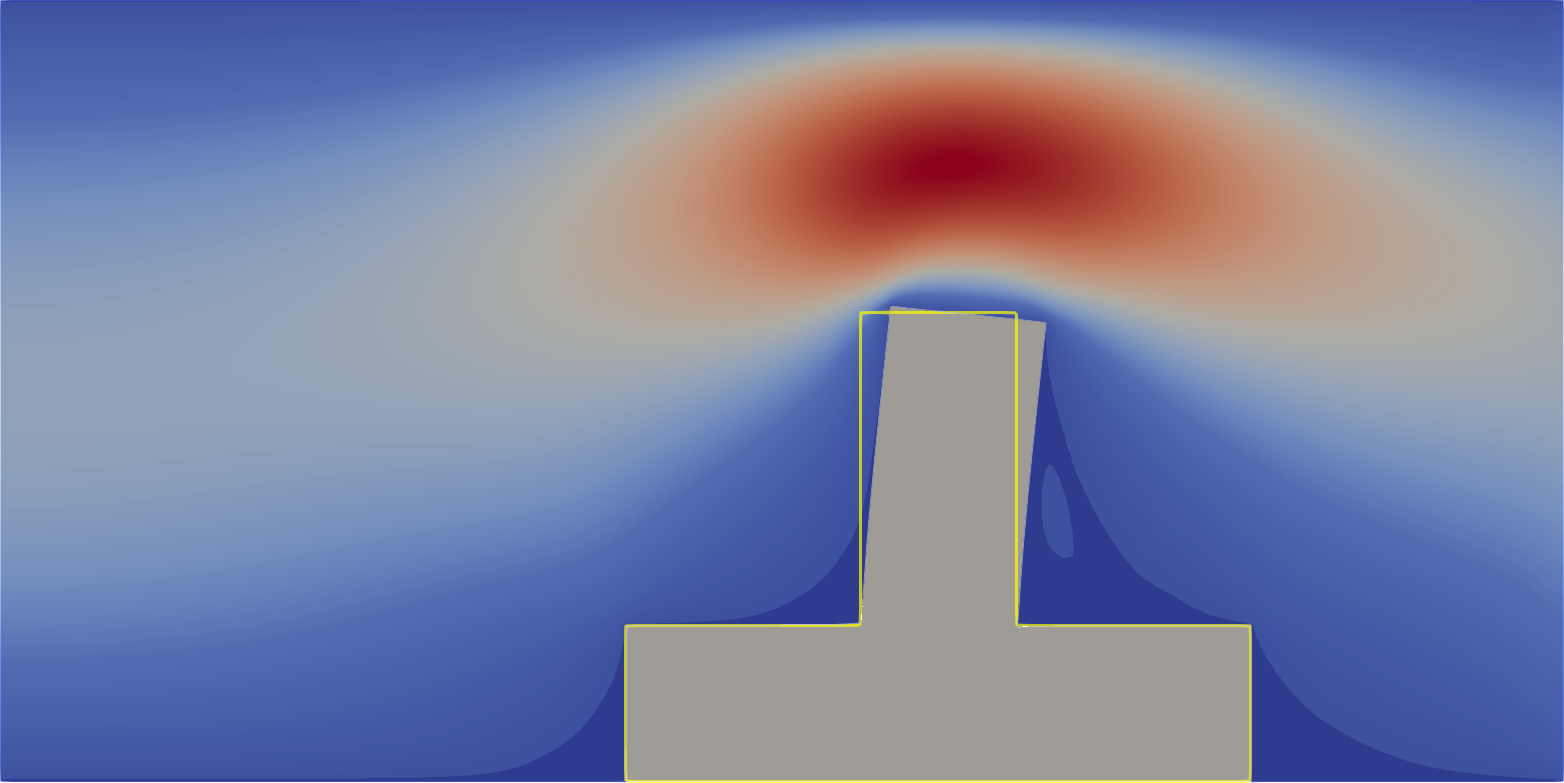}
		\caption{}
		\label{fig:picFluidStructureInitial}
	\end{subfigure}
	\begin{subfigure}{0.49\textwidth}
		\cropped{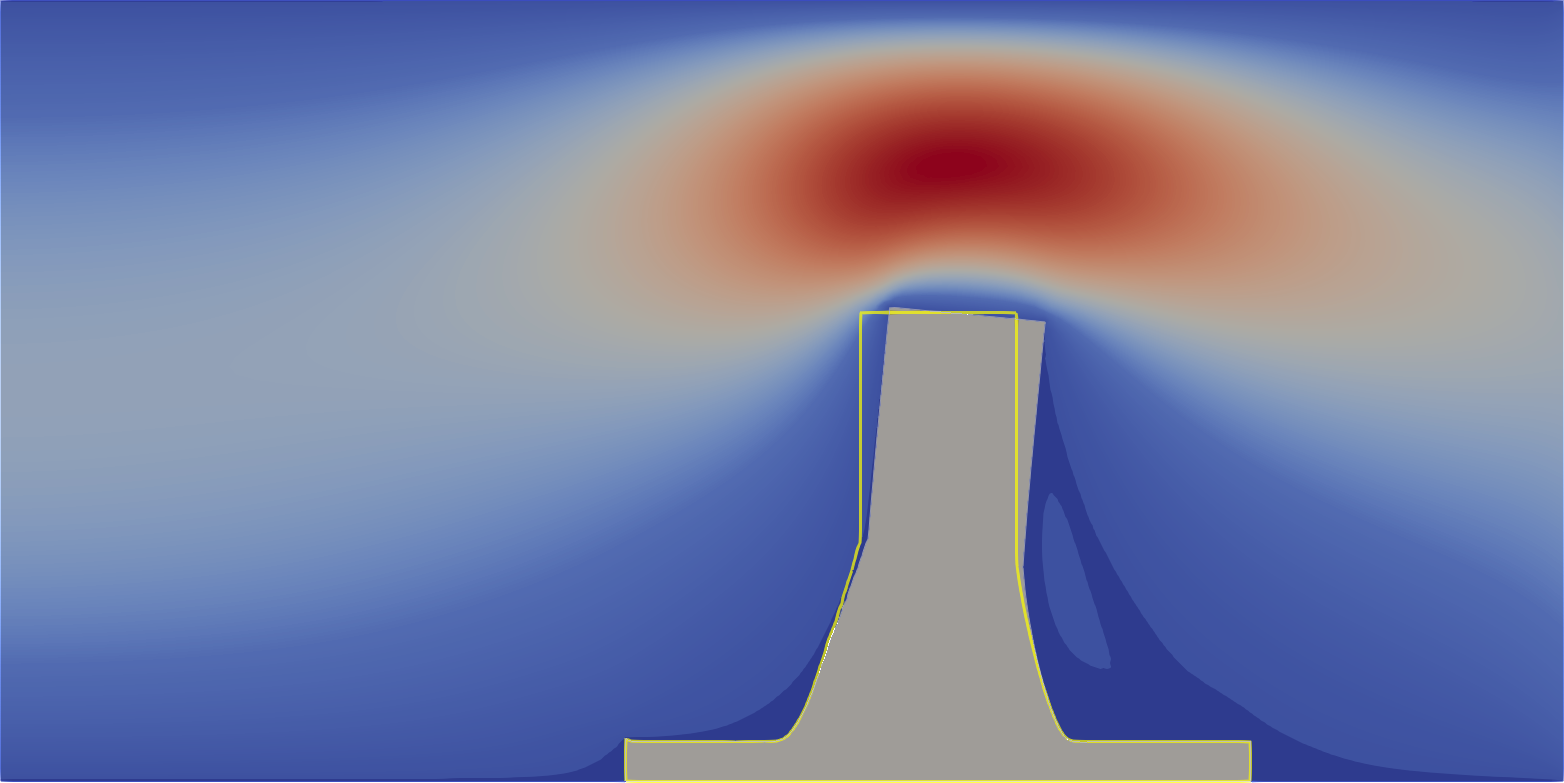}
		\caption{}
		\label{fig:picFluidStructureMin}
	\end{subfigure}
\caption{Fluid-structure interaction problem. The deformed structure is shown in grey. The undeformed structure is indicated by the yellow outlines. The colors in the fluid domain represent the norm of the fluid velocity. Red corresponds to high velocity, blue corresponds to low velocity: (\subref{fig:picFluidStructureInitial}) initialization; (\subref{fig:picFluidStructureMin}) optimized design}
\label{fig:picFluidStructure}
\end{figure*}
The initialization and the optimized material distribution are depicted in \Cref{fig:picFluidStructure}.  The evolution of the objective function is reported in \Cref{fig:convergenceFluidStructure}.

\section{Conclusion}\label{sec::conclusion}
We developed the \python package \autofreefem designed for the automatic generation of simulation code and corresponding problem documentation to facilitate the simulation and optimization of complex non-linear multi-physics problems. A \latex component enables users to produce consistent documentations, while the \freefem component focuses on the numerical simulation aspect, providing a robust platform for solving partial differential equations. The effectiveness of our approach has been demonstrated through its application to various shape optimization problems. 

We believe that this integrated approach offers several significant pedagogical advantages. Firstly, it minimizes the risk of discrepancies between the documented theory and the implemented code, as both are derived from the same underlying source. This consistency is crucial for the reproducibility of scientific results. Secondly, the automation of code and documentation generation saves time and reduces the potential for human error, especially for students or beginners in the field.

\backmatter

\bmhead{Supplementary information}

%
%

\section*{Compliance with ethical standards}

\paragraph{Funding} The authors did not receive support from any organization for the submitted work.

\paragraph{Conﬂict of interest} The authors declare that they have no conflict of interest.

\paragraph{Replication of results} The developed python package is available at \href{https://gitlab.tugraz.at/autofreefem/autofreefem}{https://gitlab.tugraz.at/autofreefem/autofreefem}. This allows to reproduce all results of the present paper. All computations were performed using \freefem version 4.14.

%
%

\begin{appendices}
\clearpage\onecolumn
	\section{\latex documentation of the linear elasticity cantilever problem}\label{appendix::Cantilever}
 All boxed content in the appendix has been automatically generated by \autofreefem. In order to demonstrate the capabilities of the software, no manual improvements have been made to the output.
 \begin{mdframed}
	
Let $\phi(x)$ be the level-set function and
 \begin{align*}\raisebox{\depth}{$\chi$}(x) &= \begin{cases} 1 & \mbox{if}\quad \phi(x) < 0 \\ 1/100 & \mbox{if}\quad \phi(x)  \ge 0 \end{cases}.\end{align*}
The Lagrangian of the linear Elasticity problem is
\begin{dmath*}
\mathcal L(\left[  {\mathbf{ u } }\right],\left[ \delta{\mathbf{ u}}\right]) = J( {\mathbf{ u } }) -  \int_{\Gamma_N} f \left( \mathbf e_y \cdot \delta{\mathbf{ u}} \right) \,ds_x +  \int_{D} \left( \pmb \sigma( {\mathbf{ u } }) \mathrel{:}  \nabla \delta{\mathbf{ u}}  \right) \,dx , 
\end{dmath*}
with 
\begin{align*}
\lambda &= - \frac{E \nu}{2 \nu^{2} + \nu - 1},\\
J( {\mathbf{ u } }) &=  \int_{D} \ell \raisebox{\depth}{$\chi$} \,dx +  \int_{D} \left( \pmb \sigma( {\mathbf{ u } }) \mathrel{:} \pmb \varepsilon( {\mathbf{ u } }) \right) \,dx,\\
\pmb \varepsilon( {\mathbf{ u } }) &= \frac{ \nabla  {\mathbf{ u } }  +   \nabla  {\mathbf{ u } } ^\top }{2},\\
\mu &= \frac{E}{2 \left(\nu + 1\right)},\\
\pmb \sigma( {\mathbf{ u } }) &= \left(\lambda \mathbf{I}  \,\operatorname{tr}\left(\pmb \varepsilon( {\mathbf{ u } })\right)  + 2 \mu \pmb \varepsilon( {\mathbf{ u } })\right) \raisebox{\depth}{$\chi$}.
 \end{align*}
The state $\left[  {\mathbf{ u } }\right]$ is the solution of the linear problem
\begin{dmath*}
-  \int_{\Gamma_N} f \left( \mathbf e_y \cdot \delta{\mathbf{ u}} \right) \,ds_x +  \int_{D} \left( \pmb \sigma( {\mathbf{ u } }) \mathrel{:}  \nabla \delta{\mathbf{ u}}  \right) \,dx = 0\qquad \forall\; \delta{\mathbf{ u}} . 
\end{dmath*}
The adjoint state $\left[ \tilde {\mathbf{  u } }\right]$ to the direct state $\left[  {\mathbf{ u } }\right]$ is the solution of
\begin{dmath*}
\partial_{( {\mathbf{ u } },\delta{\mathbf{ u}})} J( {\mathbf{ u } },\delta{\mathbf{ u}}) +  \int_{D} \left( \partial_{( {\mathbf{ u } },\delta{\mathbf{ u}})} \pmb \sigma(\delta{\mathbf{ u}}) \mathrel{:}  \nabla \tilde {\mathbf{  u } }  \right) \,dx = 0\qquad \forall\; \delta{\mathbf{ u}} , 
\end{dmath*}
with 
\begin{align*}
\partial_{( {\mathbf{ u } },\delta{\mathbf{ u}})} J( {\mathbf{ u } },\delta{\mathbf{ u}}) &=  \int_{D} \left( \partial_{( {\mathbf{ u } },\delta{\mathbf{ u}})} \pmb \sigma(\delta{\mathbf{ u}}) \mathrel{:} \pmb \varepsilon( {\mathbf{ u } }) \right) + \left( \pmb \sigma( {\mathbf{ u } }) \mathrel{:} \partial_{( {\mathbf{ u } },\delta{\mathbf{ u}})} \pmb \varepsilon(\delta{\mathbf{ u}}) \right) \,dx,\\
\partial_{( {\mathbf{ u } },\delta{\mathbf{ u}})} \pmb \sigma(\delta{\mathbf{ u}}) &= \left(\lambda \mathbf{I}  \,\operatorname{tr}\left(\partial_{( {\mathbf{ u } },\delta{\mathbf{ u}})} \pmb \varepsilon(\delta{\mathbf{ u}})\right)  + 2 \mu \partial_{( {\mathbf{ u } },\delta{\mathbf{ u}})} \pmb \varepsilon(\delta{\mathbf{ u}})\right) \raisebox{\depth}{$\chi$},\\
\partial_{( {\mathbf{ u } },\delta{\mathbf{ u}})} \pmb \varepsilon(\delta{\mathbf{ u}}) &= \frac{ \nabla \delta{\mathbf{ u}}  +   \nabla \delta{\mathbf{ u}} ^\top }{2}.
 \end{align*}
In order to compute the shape derivative, we consider a shape pertubation $\mathbf x = T_t(\mathbf X) = \mathbf X + t \mathbf  V(X)$ with a suitable velocity field $\mathbf V$. The perturbed Lagrangian $\mathcal G$ using the Lagrangian state $\left[  {\mathbf{ u } }\right]^{L}$ is given by
\begin{dmath*}
\mathcal G(t,\left[  {\mathbf{ u } }\right]^{L},\left[ \delta{\mathbf{ u}}\right]^{L}) = \mathcal V(J) -  \int_{\Gamma_N} f  \operatorname{det}\mathbf F( t )  \left( \mathbf e_y \cdot \delta{\mathbf{ u}} \right) \sqrt{\left( \left( \mathbf n \cdot  \mathbf F^{-1}(t) \right) \cdot \left( \mathbf n \cdot  \mathbf F^{-1}(t) \right) \right)} \,ds_x +  \int_{D}  \operatorname{det}\mathbf F( t )  \left( \mathcal V(\pmb \sigma) \mathrel{:} \left( \left(  \nabla \delta{\mathbf{ u}}  \cdot  \mathbf F^{-1}(t) \right) \right)\right) \,dx , 
\end{dmath*}
with
\begin{dmath*}
\mathcal V(\pmb \sigma) = \left(\lambda \mathbf{I}  \,\operatorname{tr}\left(\mathcal V(\pmb \varepsilon)\right)  + 2 \mu \mathcal V(\pmb \varepsilon)\right) \raisebox{\depth}{$\chi$},
 \end{dmath*}
\begin{dmath*}
\mathcal V(\pmb \varepsilon) = \frac{\left(  \nabla  {\mathbf{ u } }  \cdot  \mathbf F^{-1}(t) \right) +  \left( \left(  \nabla  {\mathbf{ u } }  \cdot  \mathbf F^{-1}(t) \right) \right)^\top }{2},
 \end{dmath*}
\begin{dmath*}
\mathcal V(J) =  \int_{D}  \operatorname{det}\mathbf F( t )  \left( \mathcal V(\pmb \sigma) \mathrel{:} \mathcal V(\pmb \varepsilon) \right) \,dx +  \int_{D} \ell \raisebox{\depth}{$\chi$}  \operatorname{det}\mathbf F( t )  \,dx.
 \end{dmath*}
For the direct state $\left[  {\mathbf{ u } }\right]$ and the adjoint state $\left[ \tilde {\mathbf{  u } }\right]$, the volume expression of the shape derivative is given by
\begin{dmath*}
DJ(\Omega)(\mathbf V) = \partial_{t} J( {\mathbf{ u } },{\mathbf V}) -  \int_{\Gamma_N} f  \operatorname{div} {\mathbf V}  \left( \mathbf e_y \cdot \tilde {\mathbf{  u } } \right) \,ds_x +  \int_{\Gamma_N} f \left( \mathbf n \cdot \left( \mathbf n \cdot  \nabla {\mathbf V}  \right) \right) \left( \mathbf e_y \cdot \tilde {\mathbf{  u } } \right) \,ds_x +  \int_{D}  \operatorname{div} {\mathbf V}  \left( \pmb \sigma( {\mathbf{ u } }) \mathrel{:}  \nabla \tilde {\mathbf{  u } }  \right) + \left( \partial_{t} \pmb \sigma( {\mathbf{ u } },{\mathbf V}) \mathrel{:}  \nabla \tilde {\mathbf{  u } }  \right) - \left( \pmb \sigma( {\mathbf{ u } }) \mathrel{:} \left( \left(  \nabla \tilde {\mathbf{  u } }  \cdot  \nabla {\mathbf V}  \right) \right)\right) \,dx , 
\end{dmath*}
with
\begin{dmath*}
\partial_{t} J( {\mathbf{ u } },{\mathbf V}) =  \int_{D} \ell \raisebox{\depth}{$\chi$}  \operatorname{div} {\mathbf V}  \,dx +  \int_{D}  \operatorname{div} {\mathbf V}  \left( \pmb \sigma( {\mathbf{ u } }) \mathrel{:} \pmb \varepsilon( {\mathbf{ u } }) \right) + \left( \partial_{t} \pmb \sigma( {\mathbf{ u } },{\mathbf V}) \mathrel{:} \pmb \varepsilon( {\mathbf{ u } }) \right) + \left( \pmb \sigma( {\mathbf{ u } }) \mathrel{:} \partial_{t} \pmb \varepsilon( {\mathbf{ u } },{\mathbf V}) \right) \,dx,
 \end{dmath*}
\begin{dmath*}
\partial_{t} \pmb \sigma( {\mathbf{ u } },{\mathbf V}) = \left(\lambda \mathbf{I}  \,\operatorname{tr}\left(\partial_{t} \pmb \varepsilon( {\mathbf{ u } },{\mathbf V})\right)  + 2 \mu \partial_{t} \pmb \varepsilon( {\mathbf{ u } },{\mathbf V})\right) \raisebox{\depth}{$\chi$},
 \end{dmath*}
\begin{dmath*}
\partial_{t} \pmb \varepsilon( {\mathbf{ u } },{\mathbf V}) = - \frac{\left(  \nabla  {\mathbf{ u } }  \cdot  \nabla {\mathbf V}  \right) +  \left( \left(  \nabla  {\mathbf{ u } }  \cdot  \nabla {\mathbf V}  \right) \right)^\top }{2}.
 \end{dmath*}
       
 \end{mdframed}

\end{appendices}
\twocolumn
\bibliography{literature}

\end{document}